\documentclass[numbers,webpdf,imanum]{ima-authoring-template}%

\usepackage{color}
\usepackage{hyperref}
\hypersetup{colorlinks=true,linkcolor=blue!70!black,filecolor=purple!50!black,citecolor=blue!50!green}
\usepackage{url}
\usepackage{accents}
\usepackage{amsthm}
\usepackage{thmtools}
\usepackage[noabbrev,nameinlink,capitalise]{cleveref}
\usepackage{enumerate}

\newcommand{\dd}{\mathrm{d}}

\newcommand{\ol}[1]{\overline{#1}}
\newcommand{\ull}[1]{\underline{#1}}
\newcommand{\spl}{\langle}
\newcommand{\spr}{\rangle}
\newcommand{\bpm}{\begin{pmatrix}}
\newcommand{\epm}{\end{pmatrix}}
\newenvironment{equations}{\equation\aligned}{\endaligned\endequation}

\newcommand{\dom}{\mathcal{O}} %
\newcommand{\nv}{\boldsymbol{\nu}}
\newcommand{\bflow}{\mathbf{b}}

\newcommand{\angvel}{\Omega}
\newcommand{\opd}{(\omega+i\conv+i\angvel\times)}
\newcommand{\acow}{a} %
\renewcommand{\u}{\mathbf{u}}
\renewcommand{\v}{\mathbf{v}}
\newcommand{\w}{\mathbf{w}}
\newcommand{\z}{\mathbf{z}}
\newcommand{\m}{\underline{\underline{m}}} %
\newcommand{\q}{\mathbf{q}} %
\newcommand{\M}{M}
\newcommand{\F}{F}
\newcommand{\betad}{\beta_{\text{disc}}}

\newcommand{\nh}{n} %
\newcommand{\limnh}{n\to\infty} %
\newcommand{\Bstab}{B}
\newcommand{\Bcomp}{K}

\newcommand{\Hz}{\bH^1_{\nv0}}
\newcommand{\Hzz}{\bH^1_0}
\newcommand{\X}{\IX}
\newcommand{\Xz}{\X} %
\newcommand{\Xn}{\X_n}
\newcommand{\Xnz}{\Xn} %
\newcommand{\Qn}{Q_n}
\newcommand{\Qnp}{Q_n^+}

\newcommand{\Vnn}{V}
\newcommand{\Wnn}{W}
\newcommand{\TWnn}{\widetilde{W}}
\newcommand{\Vhnn}{\Vnn_n}
\newcommand{\Whnn}{\Wnn_n}
\newcommand{\TWhnn}{\widetilde \Wnn_n}

\newcommand{\D}{D} %
\newcommand{\PXnz}{P_{\Xnz}} %
\newcommand{\PQn}{P_{\Qn}} %
\newcommand{\PQnp}{P_{\Qnp}} %

\newcommand{\PV}{P_\Vnn}
\newcommand{\PW}{P_\Wnn}
\newcommand{\Tn}{T_n}
\newcommand{\dof}{\mathrm{dof}}

\newcommand{\Cps}{C_{PS}} %
\newcommand{\Cxh}{C_{XH}} %
\newcommand{\coeff}{\alpha}

\newcommand{\mesh}{{\mathcal{T}_n}} 

\renewcommand{\Re}{\operatorname{Re}}
\renewcommand{\Im}{\operatorname{Im}}
\newcommand{\conv}{\partial_{\bflow}} %
\renewcommand{\div}{\operatorname{div}}
\DeclareMathOperator{\grad}{grad}

\DeclareMathOperator{\hess}{Hess}

\DeclareMathOperator{\supp}{supp}
\DeclareMathOperator{\spn}{span}

\renewcommand{\dim}{\operatorname{dim}}
\DeclareMathOperator{\codim}{codim}
\DeclareMathOperator{\ran}{ran}

\DeclareMathOperator{\sign}{sgn}

\DeclareMathOperator{\mean}{mean}

\newcommand{\inner}[1]{\langle #1 \rangle}

\newcommand*{\norm}[1]{\left\|#1\right\|}
\newcommand{\mnull}{\setminus \{0\}}
\newcommand\restr[2]{{ \left.\kern-\nulldelimiterspace #1 \vphantom{\big|} \right|_{#2} }}

\newcommand{\IC}{\mathbb{C}}

\newcommand{\IN}{\mathbb{N}}

\newcommand{\IR}{\mathbb{R}}

\newcommand{\IX}{\mathbb{X}}

\newcommand{\calO}{\mathcal{O}}

\newcommand{\calT}{\mathcal{T}}

\newcommand{\be}{\mathbf{e}}
\newcommand{\bff}{\mathbf{f}}

\newcommand{\bt}{\mathbf{t}}
\newcommand{\bu}{\mathbf{u}}
\newcommand{\bv}{\mathbf{v}}

\newcommand{\bx}{\mathbf{x}}
\newcommand{\by}{\mathbf{y}}

\newcommand{\bB}{\mathbf{B}}
\newcommand{\bC}{\mathbf{C}}

\newcommand{\bH}{\mathbf{H}}

\newcommand{\bL}{\mathbf{L}}
\newcommand{\bM}{\mathbf{M}}

\newcommand{\bS}{\mathbf{S}}

\newcommand{\bX}{\mathbf{X}}

\definecolor{pscol}{rgb}{0.8,0,0}

\definecolor{mhcol}{rgb}{0.0,0.6,0}

\definecolor{revcol}{rgb}{0.7,0.4,0}

\usepackage{stmaryrd} %
\usepackage{pgf,tikz}
\usetikzlibrary{spy,matrix, calc, arrows,shapes,decorations}
\usepackage{tkz-euclide}
\usetikzlibrary{positioning,arrows,calc,decorations.markings,math,arrows.meta}
\usepackage{pgfplots} 
\usepgfplotslibrary{groupplots}
\usepackage{csvsimple} 
\usepackage{siunitx} %
\usepackage{tabularx}
\usepackage{booktabs}
\usepackage{multirow}
\usepackage[normalem]{ulem}
\useunder{\uline}{\ul}{}

\pgfplotscreateplotcyclelist{paulcolors}{%
violet,every mark/.append style={solid,fill=violet},mark=square*,very thick,mark size=3pt\\%
teal,every mark/.append style={solid,fill=teal},mark=*,very thick,mark size=3pt\\%
orange,every mark/.append style={solid,fill=orange},mark=diamond*,very thick,mark size=3pt\\%
cyan,every mark/.append style={solid,fill=cyan},mark=star,very thick,mark size=3pt\\%
}
\pgfplotscreateplotcyclelist{paulcolors2}{%
teal,every mark/.append style={solid,fill=teal},mark=*,very thick,mark size=3pt\\%
orange,every mark/.append style={solid,fill=orange},mark=diamond*,very thick,mark size=3pt\\%
teal,densely dashed,every mark/.append style={solid,fill=teal},mark=*,very thick,mark size=3pt\\%
orange,densely dashed,every mark/.append style={solid,fill=orange},mark=diamond*,very thick,mark size=3pt\\%
teal,dash dot,every mark/.append style={solid,fill=teal},mark=*,very thick,mark size=3pt\\%
orange,dash dot,every mark/.append style={solid,fill=orange},mark=diamond*,very thick,mark size=3pt\\%
}
\pgfplotscreateplotcyclelist{paulcolors1}{%
teal,every mark/.append style={solid,fill=teal},mark=*,very thick,mark size=3pt\\%
}
\pgfplotscreateplotcyclelist{paulcolorsfill}{%
teal,fill=teal,fill opacity=0.5,thick\\
orange,fill=orange,fill opacity=0.5,thick\\
violet,fill=violet,fill opacity=0.5,thick\\
}

\pgfplotsset{
    discard if not/.style 2 args={
        x filter/.append code={
            \edef\tempa{\thisrow{#1}}
            \edef\tempb{#2}
            \ifx\tempa\tempb
            \else
                
            \fi
        }
    }
}
\pgfplotsset{compat=1.16}%

\pgfplotsset{ width=.42\linewidth}

\newcommand{\logLogSlopeTriangle}[5]
{

    \pgfplotsextra
    {
        \pgfkeysgetvalue{/pgfplots/xmin}{\xmin}
        \pgfkeysgetvalue{/pgfplots/xmax}{\xmax}
        \pgfkeysgetvalue{/pgfplots/ymin}{\ymin}
        \pgfkeysgetvalue{/pgfplots/ymax}{\ymax}

        \pgfmathsetmacro{\xArel}{#1}
        \pgfmathsetmacro{\yArel}{#3}
        \pgfmathsetmacro{\xBrel}{#1-#2}
        \pgfmathsetmacro{\yBrel}{\yArel}
        \pgfmathsetmacro{\xCrel}{\xArel}

        \pgfmathsetmacro{\lnxB}{\xmin*(1-(#1-#2))+\xmax*(#1-#2)} %
        \pgfmathsetmacro{\lnxA}{\xmin*(1-#1)+\xmax*#1} %
        \pgfmathsetmacro{\lnyA}{\ymin*(1-#3)+\ymax*#3} %
        \pgfmathsetmacro{\lnyC}{\lnyA-#4*(\lnxA-\lnxB)}
        \pgfmathsetmacro{\yCrel}{\lnyC-\ymin)/(\ymax-\ymin)} %

        \coordinate (A) at (rel axis cs:\xArel,\yArel);
        \coordinate (B) at (rel axis cs:\xBrel,\yBrel);
        \coordinate (C) at (rel axis cs:\xCrel,\yCrel);

        \draw[#5]   (A)-- node[pos=0.5,anchor=north] {}
                    (B)-- 
                    (C)-- node[pos=0.5,anchor=west] {#4}
                    cycle;
    }
}

\newcommand\convergence{\hyperlink{converge}{convergence${}^{i)}$}}
\newcommand\stable{\hyperlink{stable}{stable${}^{v)}$}}
\newcommand\regular{\hyperlink{regular}{regular${}^{vi)}$}}
\newcommand\regularity{\hyperlink{regular}{regularity${}^{vi)}$}}

\newcommand\compactop{\hyperlink{compactop}{compact${}^{iv)}$}}
\newcommand\compactnessseq{\hyperlink{compactseq}{compactness${}^{ii)}$}}

\newcommand\approximates{\hyperlink{approximate}{approximates${}^{iii)}$}}
\newcommand\asymptoticconsistency{\hyperlink{approximate}{asymptotic consistency${}^{iii)}$}}

\usepackage{comment}
\excludecomment{extproof}

\theoremstyle{thmstyletwo}%
\newtheorem{theorem}{Theorem}%
\newtheorem{cor}[theorem]{Corollary}
\newtheorem{lemma}[theorem]{Lemma}
\newtheorem{proposition}[theorem]{Proposition}%
\newtheorem{assumption}[theorem]{Assumption}

\newtheorem{remark}{Remark}%

\newtheorem{definition}{Definition}

\numberwithin{equation}{section}

\begin{document}

\DOI{DOI HERE}
\copyrightyear{2021}
\vol{00}
\pubyear{2021}
\access{Advance Access Publication Date: Day Month Year}
\appnotes{Paper}
\copyrightstatement{Published by Oxford University Press on behalf of the Institute of Mathematics and its Applications. All rights reserved.}
\firstpage{1}

\title[A new T-compatibility condition applied to Galbrun's equation]
{A new T-compatibility condition and its application to the discretization of the damped time-harmonic Galbrun's equation}

\author{Martin Halla
\address{
    \orgdiv{Institut f\"ur Angewandte und Numerische Mathematik},
	\orgname{Karlsruher Institut f\"ur Technologie},
	\orgaddress{\street{Englerstr.\ 2}, \postcode{76131} \state{Karlsruhe}, \country{Deutschland}} 
}}

\author{Christoph Lehrenfeld
\address{
    \orgdiv{Institut f\"ur Numerische und Angewandte Mathematik},
    \orgname{Georg-August Universität Göttingen},
    \orgaddress{\street{Lotzestr.\ 16-18}, \postcode{37083} \state{Göttingen}, \country{Deutschland}} 
}}

\author{Paul Stocker*
\address{
    \orgdiv{Fakult\"at für Mathematik, Universit\"at Wien},
    \orgaddress{\street{Oskar-Morgenstern-Platz 1}, \postcode{1090} \state{Wien}, \country{\"Osterreich}}
}}
\corresp[*]{Corresponding author: \href{email:paul.stocker@univie.ac.at}{paul.stocker@univie.ac.at}}

\authormark{Martin Halla, Christoph Lehrenfeld, Paul Stocker}

\received{Date}{0}{Year}
\revised{Date}{0}{Year}
\accepted{Date}{0}{Year}

\abstract{
We consider the approximation of weakly T-coercive operators.
The main property to ensure the convergence thereof is the regularity of the approximation (in the vocabulary of discrete approximation schemes).
In a previous work the existence of discrete operators $T_n$ which converge to $T$ in a discrete norm was shown to be sufficient to obtain regularity.
Although this framework proved useful for many applications for some instances the former assumption is too strong.
Thus in the present article we report a weaker
{criterion}
for which the discrete operators $T_n$ only have to converge point-wise, but in addition a weak T-coercivity condition has to be satisfied on the discrete level.
We apply the new framework to prove the convergence of certain $H^1$-conforming finite element discretizations of the damped time-harmonic Galbrun's equation, which is used to model the oscillations of stars.
A main ingredient in the latter analysis is the uniformly stable invertibility of the divergence operator on certain spaces, which is related to the topic of stable discretizations of the Stokes equation.
}

\keywords{discrete approximation schemes, weak T-coercivity, Galbrun's equation\\ 
\textit{MSC Classification:} {35L05, 35Q85, 65N30}}

\maketitle

\section{Introduction}

An origin of the T-coercivity technique to analyze equations of non weakly coercive form can be found in the theory of Maxwell's equation and goes back at least to \cite{BuffaCostabelSchwab:02,Buffa:05}.
The idea to use a discrete variant to prove the stability of approximations can be found e.g.\ in \cite{HohageNannen:15,BonnetBDCarvalhoCiarlet:18}.
In \cite{Halla:21Tcomp} this approach was formalized to a framework to prove the convergence of Galerkin approximations of holomorphic eigenvalue problems and was successfully applied for perfectly matched layer methods to scalar isotropic \cite{Halla:21PML} and anisotropic \cite{Halla:22PMLani} materials, Maxwell problems in conductive media \cite{Halla:21Tcomp}, modified Maxwell Steklov problems \cite{Halla:21SteklovAppr} and Maxwell transmission problems for dispersive media \cite{Unger:21,Halla:21SC}.
In particular \cite{Halla:21Tcomp} is build upon the much broader framework of \emph{discrete approximation schemes} \cite{Stummel:70,Vainikko:76} which originated in the 1970s and the best results for eigenvalue problems in this context are \cite{Karma:96a,Karma:96b}.
The main contribution of \cite{Halla:21Tcomp} was to provide a practical criterion to prove the regularity of approximations, which allows to apply the results achieved for discrete approximation schemes.
Although for some applications it turns out that the T-compatibility criterion of \cite{Halla:21Tcomp}  is too strong, and hence we present in this article a weaker variant.
Some similarity can be drawn to the analysis of p-finite element methods for Maxwell problems \cite{BCDDH:11}, for which (opposed to h-finite element methods) the cochain projections are not uniformly $L^2$ bounded, and hence the discrete compactness property is obtained in \cite{BCDDH:11} by an alternative technique.

{Primarily the T-coercivity approach serves a technique for the analysis of PDEs and the numerical analysis of respective discretizations.
However, $T$-coercivity techniques can also be used to construct new numerical schemes. Indeed, if feasible, the operator $T$ can be included in the discretized variational formulation as e.g.\ done in \cite{CiarletJamelot24,HallaHohageOberender24}. Having now the discretization of a weakly coercive problem at hand, the stability of the approximations follow in a straightforward manner.}

The present article is motivated by the study of approximations to the damped time-harmonic Galbrun's equation.
The Galbrun's equation~\cite{Galbrun:31} is a linearization of the nonlinear Euler equations with the Lagrangian perturbation of displacement as unknown, and is used in aeroacoustics \cite{maeder202090} as well as in an extended version in asteroseismology \cite{LyndenBOstriker:67}. %
We refer to \cite{HaeggBerggren:19} for a well-posedness analysis in the time domain.
In the time-harmonic domain an approach in aeroacoustics is to use a stabilized formulation, which is justified by the introduction of an additional transport equation for the vorticity, and we refer to the well-posedness analysis in \cite{BonnetBDMercierMillotPernetPeynaud:12}.
Different to aeroacoustics in asteroseismology there exists a significant damping of waves which allows the equation to be analyzed in a more direct way, see the well-posedness results \cite{HallaHohage:21,Halla:21GalExt}.
In the second part of the present article we apply our new framework to the approximation of the damped time-harmonic Galbrun's equation as considered in \cite{HallaHohage:21}:
\begin{align}\label{eq:Galbrun}
\begin{aligned}
&-\rho\opd^2\u 
- \nabla\left(\rho c_s^2\div \u\right)
+ (\div \u) \nabla p\\
&\qquad\qquad-\nabla(\nabla p\cdot \u) + (\hess(p)-\rho\hess(\phi))\u 
+ \gamma \rho (-i \omega) \u
= \bff \quad \mbox{in } \dom
\end{aligned}
\end{align}
where $\rho, p, \phi, c_s,\bflow, \angvel$ and $\bff$ denote density, pressure, gravitational potential, sound speed, background velocity, angular velocity of the frame and sources, $\partial_\bflow := \sum_{l=1}^3 \bflow_l\partial_{x_l}$ denotes the directional derivative in direction $\bflow$, $\hess(p)$ the Hessian of $p$, $\dom\subset\mathbb{R}^3$ a bounded domain, and damping is modeled by the term $- i \omega \gamma \rho \u$ with damping coefficient $\gamma$.
The main challenge to tackle this equation can already be observed in the case $p, \phi=0, \angvel=0$.
We discretize \eqref{eq:Galbrun} with conforming $\bH^1$ finite elements.
To guarantee the stability of the approximation we use vectorial finite element spaces which admit a suitable uniformly stable inversion of the divergence operator.
In particular, let $X_h\subset\bH^1$ be a Lagrangian vectorial finite element space of order $k$ and $Q_h\subset L^2$ be a scalar finite element space 
and $L^2_0:=\{u\in L^2\colon \int_\Omega u\,\dd x=0\}$.
Then we require that there exists a uniformly bounded inverse of the (discrete) divergence operator acting on the spaces $Q_h\cap L^2_0\to X_h\cap\bH^1_0$.
Such methods have been developed in the field of computational fluid dynamics 
for the stable discretization of incompressible Stokes and Navier-Stokes equations, cf. e.g. \cite{John:16}.

Especially convenient for the analysis are the so-called \emph{divergence free finite elements}, meaning that the approximative solutions to the Stokes equations are exactly divergence free.
However, note that there exist sophisticated techniques to construct such elements and not all \emph{divergence free finite elements} fit our needs. %
The pioneering work for \emph{divergence free finite elements} was set by Scott and Vogelius~\cite{ScottVogelius:85}, who established respective results (suitable for our purpose) in 2D for triangular quasi-uniform meshes with \emph{finite degree of degeneracy} and polynomial degree $k\geq 4$ (the quasi-uniformity is actually not necessary due to \cite{FalkNeilan:13}).
In three dimensions Zhang \cite{Zhang:11a} reported a generalization to uniform tetrahedral grids for $k\geq6$, and his results in \cite{Zhang:09} indicate that for general tetrahedral grids suitable orders are $k\geq8$.
The application of convenient finite element spaces on specialized meshes generated by barycentric refinements (suitable for our purpose)
received extensive attention and we refer e.g.\ to \cite{ArnoldQin:92,Zhang:05,GuzmanNeilan:18}. %
In general such schemes are related to respective discretizations of suitable deRahm complexes with high regularity \cite{ChristiansenHu:18,Neilan:20,GuzmanLischkeNeilan:20}.
There exist also several results for elements on quadrilateral grids for which we refer to the bibliographies of \cite{JohnLinkeMerdonNeilanRebholz:17,Neilan:20}.
Other approaches to construct \emph{divergence free finite elements} include enriched finite elements, nonconforming elements, discontinuous Galerkin methods and isogeometric methods.

Although we will make use of the advantages of \emph{divergence free finite elements} in the analysis, we note that the more important property is the stable Stokes approximation.
A comparison and analysis of different robust finite element discretizations for a simplified Galbrun's equation is presented in \cite{2205.15650}.
{Approximations with $H(\div)$-conforming finite elements and DGFEMs  are analyzed in \cite{2306.03496,Thesis_vB23,HallaLehrenfeldvanBeeck25}, employing the framework of the current article}.

The remainder of this article is structured as follows.
In \Cref{sec:framework} we report a multipurpose framework based on a weak T-compatibility condition (weaker than in \cite{Halla:21Tcomp}) to obtain the regularity and hence the stability of approximations.
Although in the present article we consider only conforming discretizations to \eqref{eq:Galbrun}, we formulate the framework in a general way to include also nonconforming approximations {\cite{2306.03496}}.
In \Cref{sec:galbrun} we apply the former framework to discretizations of \eqref{eq:Galbrun}.
In particular, in \Cref{subsec:homo} we consider a simplified case of \eqref{eq:Galbrun} to present the main ideas and in \Cref{subsec:hetero} we treat the general case.
In \Cref{sec:num} we present computational examples to accompany our theoretical results and we conclude in \Cref{sec:conclusion}.

\section{Abstract framework}\label{sec:framework}
This section discusses a multipurpose framework for the analysis of approximations of linear operators. 
In \Cref{subsec:DAS} we review the framework and important definitions, as well as sufficient conditions for the convergence of the approximative solution.
We aim to apply this framework to operators that are Fredholm with index zero, however have the structure of `coercive+compact' only up to a bijection. 
Such operators are called weakly $T$-coercive, a precise definition is given in \Cref{subsec:newTcomp}.
Note that this property is equivalent to an operator being Fredholm with index zero, and the construction of a suitable $T$ operator is the tool to prove this property.
Here we study a way how this property can be mimiced on the discrete level to ensure convergent approximations.

\subsection{Discrete approximation schemes}\label{subsec:DAS}

We consider discrete approximation schemes in Hilbert spaces.
Note that the forthcoming setting is a bit more restrictive than the schemes considered in \cite{Stummel:70,Vainikko:76,Karma:96a}, but more convenient for our purposes.
For two Hilbert spaces $(X, \inner{\cdot,\cdot}_X)$, $(Y, \inner{\cdot,\cdot}_Y)$ let $L(X,Y)$ be the space of bounded linear operators from $X$ to $Y$, and set $L(X):=L(X,X)$.

\begin{definition}\label{def:DAS}
We call $\{X_n, A_n, p_n\}_{n\in\IN}$ a \emph{discrete approximation scheme} of $A\in L(X)$ if the following properties hold:
Let $(X_n, \inner{\cdot,\cdot}_{X_n})_{n\in\IN}$ be a sequence of finite dimensional Hilbert spaces and $A_n\in L(X_n)$.
And let $p_n \in L(X,X_n)$ such that $\lim_{n\to\infty} \|p_n u\|_{X_n}=\|u\|_X$ for each $u\in X$.
We then define the following properties of a discrete approximation scheme:
\begin{enumerate}[i)\hspace{0.5em}]%
\item \hypertarget{converge}A sequence $(u_n)_{n\in\IN}, u_n\in X_n$ is said to \emph{converge} to $u\in X$, if $\lim_{n\to\infty}\|p_nu-u_n\|_{X_n}=0$.
\item \hypertarget{compactseq}A sequence $(u_n)_{n\in\IN}, u_n\in X_n$ is said to be \emph{compact}, if for every subsequence $\IN'\subset\IN$ exists a subsubsequence $\IN''\subset\IN'$ such that $(u_n)_{n\in\IN''}$ converges (to a $u\in X$).
\item \hypertarget{approximate}A sequence of operators $(A_n)_{n\in\IN}, A_n\in L(X_n)$ is said to \emph{approximate} $A\in L(X)$, if $\lim_{n\to\infty}\|A_np_nu-p_nAu\|_{X_n}=0$.
In a finite element vocabulary it might be more convenient to denote this property as \emph{asymptotic consistency}.
\item \hypertarget{compactop}A sequence of operators $(A_n)_{n\in\IN}, A_n\in L(X_n)$ is said to be \emph{compact}, if for every bounded sequence $(u_n)_{n\in\IN}, u_n\in X_n$, $\|u_n\|_{X_n}\leq C$ the sequence $(A_nu_n)_{n\in\IN}$ is compact.
\item \hypertarget{stable}A sequence of operators $(A_\nh)_{\nh\in\IN}, A_\nh\in L(X_\nh)$ is said to be \emph{stable}, if there
 exist constants $C,\nh_0>0$
	such that $A_\nh$ is invertible and $\|A_\nh^{-1}\|_{L(X_\nh)} \leq C$ for all
	$\nh>\nh_0$.
\item \hypertarget{regular}A sequence of operators $(A_n)_{n\in\IN}, A_n\in L(X_n)$ is said to be \emph{regular}, if $\|u_n\|_{X_n}\leq C$ and the compactness of $(A_nu_n)_{n\in\IN}$ implies the compactness of $(u_n)_{n\in\IN}$.
\end{enumerate}
\end{definition}
Note that we do not demand that the spaces $X_n$ are subspaces of $X$. 
Instead we demand the existence of the projection operators $p_n$.
{The vocabulary introduced in \cref{def:DAS} and used throughout the manuscript may not be familiar to every reader. We hence refer to the corresponding properties of a discrete approximation scheme with an upper index linking to the corresponding property in \cref{def:DAS}.}
The central properties we are looking for in a discrete approximation scheme are \regularity~and \asymptoticconsistency~ which are sufficient for the \convergence~of  discrete solutions.
To emphasize this we recall in the following some well known results.
\begin{lemma}\label{lem:stable}
Let $A\in L(X)$ be bijective and $(A_\nh)_{\nh\in\IN}$, $A_\nh\in L(X_\nh)$ be a discrete approximation scheme which is \regular~and \approximates~$A$.
Then $(A_\nh)_{\nh\in\IN}$ is \stable.
\end{lemma}
\begin{proof}
Follows from statement 3) of \cite[Theorem~2]{Karma:96a}.
{We give a proof for the sake of completeness. 
Assume the contrary, i.e. there exists a normalized sequence $u_n\in X_n$ with $A_nu_n=0$.
Since $A_n$ is regular, there exists $u\in X$ and subsequence which we do not rename, such that $\lim\|u_n-p_nu\|=0$.
As $A_n$ approximates $A$, we have $\lim\|A_np_nu-p_nAu\|_{X_n}=0$.
Therefore, $\lim\|p_nAu\|_{X_n}=0$. 
Since $A$ is injective, it follows $u=0$.
This contradicts $\|u_n\|_{X_n}=1$, $\lim_{n\to\infty} \|u_n-p_nu\|_{X_n}=0$, and hence the claim is proven.
}
\end{proof}
\begin{lemma}\label{lem:conv}
Let $A\in L(X)$ be bijective and $(A_\nh)_{\nh\in\IN}$, $A_\nh\in L(X_\nh)$ be a discrete approximation scheme which is \stable~ and \approximates~$A$.
{Let $n_0>0$ be such that $A_\nh$ is bijective for all $n>n_0$,
$u,u_n$ be the solutions to $Au=f$ and $A_n u_n=f_n\in X_n$, and {assume that} $\lim_{n\to\infty} \|p_nf-f_n\|_{X_n}=0$.}
Then $\lim_{n\to\infty} \|p_nu-u_n\|_{X_n}=0$.
If the approximation is a conforming Galerkin scheme, i.e.\ $X_n\subset X$ and $p_n$ is the orthogonal projection onto $X_n$, $f_n=p_nf$,
{then there exists a constants $C>0$ such that}
$\|u-u_n\|_{X}\leq C \inf_{u_n'\in X_n} \|u-u_n'\|_X$
for all $n>n_0$.
\end{lemma}
\begin{proof}
{Using that $A_\nh$ is stable, i.e. has a bounded inverse, followed by the triangle inequality,
we estimate}
\begin{align*}
\|p_nu-u_n\|_{X_n}
&\leq \sup_{m>n_0} \|A_m^{-1}\|_{L(X_m)}  \|A_np_nu-A_nu_n\|_{X_n}\\
&\leq \sup_{m>n_0} \|A_m^{-1}\|_{L(X_m)} \big( \|A_np_nu-p_nAu\|_{X_n}+\|p_nAu-A_nu_n\|_{X_n} \big) \\
&= \sup_{m>n_0} \|A_m^{-1}\|_{L(X_m)} \big( \|A_np_nu-p_nAu\|_{X_n}+\|p_nf-f_n\|_{X_n} \big).
\end{align*}
It holds that $\lim_{n\to\infty}\|A_np_nu-p_nAu\|_{X_n}=0$, because $(A_\nh)_{\nh\in\IN}$ \approximates~$A$, and that $\lim_{n\to\infty}\|p_nf-f_n\|_{X_n}=0$ by assumption.
Hence the first claim is proven.

{For the second claim we recall that we are in the setting of a conforming Galerkin scheme.
We estimate using the triangle inequality, the stability of $A_\nh$ and the definition of the projection $p_n$, to obtain}
\begin{align*}
\|u-u_n\|_{X} &\leq \|u-p_nu\|_{X}+\|p_nu-u_n\|_{X}\\
&\leq \|u-p_nu\|_{X}+ \sup_{m>n_0} \|A_m^{-1}\|_{L(X_m)} \|p_nAp_nu-p_nAu\|_{X}\\
&\leq \|u-p_nu\|_{X}+ \sup_{m>n_0} \|A_m^{-1}\|_{L(X_m)} \|A\|_{L(X)} \|u-p_nu\|_{X}\\
&= \big(1+\sup_{m>n_0} \|A_m^{-1}\|_{L(X_m)} \|A\|_{L(X)}\big)  \|u-p_nu\|_{X}.
\end{align*}
\end{proof}

\subsection{The new T-compatibility condition}\label{subsec:newTcomp}
\begin{definition}\label{def:weakTcoerc}
We define the following properties for an operator $A$.
\begin{enumerate}[i)\hspace{0.5em}]
\item An operator $A\in L(X)$ is called coercive, if there exists a constant $C>0$ such that $|\inner{Au,u}_X|\geq C\|u\|_X^2$ for all $u\in X$.
\item An operator $A\in L(X)$ is called weakly coercive, if there exists a compact operator $K\in L(X)$ such that $A+K$ is coercive.
\item An operator $A$ is called (weakly) right $T$-coercive, if $T\in L(X)$ is bijective and $AT$ is (weakly) coercive.
\end{enumerate}
\end{definition}
{Our definition of weak $T$-coercivity is in spirit equivalent to the generalized Gårding inequality in \cite[Prop.\ 3]{BuffaCostabelSchwab:02}.
The generalized Gårding inequality in [11] follows from our definition of weak T-coercivity by applying the triangle inequality.
However, the reverse direction seems to require an additional argument.}

The next theorem provides a sufficient setting for a discrete approximation of a (weakly) right $T$-coercive operator. This theorem is key for the discretization and its analysis in \cref{sec:galbrun}.
\begin{theorem}\label{thm:Tcomp}
{Let sequences $(A_\nh)_{\nh\in\IN},$ $(T_\nh)_{\nh\in\IN},$ $(B_\nh)_{\nh\in\IN},$ $(K_\nh)_{\nh\in\IN}$ and $B,T\in L(X)$ 
	satisfy the following:
There exists a constant $C>0$ such that
}
for each $\nh\in\IN$
it holds $A_\nh,T_\nh,B_\nh,K_\nh\in L(X_\nh)$,
$\|T_\nh\|_{L(X_\nh)}, \|T_\nh^{-1}\|_{L(X_\nh)}, \|\Bstab_\nh\|_{L(X_\nh)}, \|\Bstab_\nh^{-1}\|_{L(X_\nh)} \leq C$, $B$ is bijective, $(K_\nh)_{\nh\in\IN}$ is \compactop~and
\begin{subequations}
  \begin{align}
\label{eq:thm3ab}
  \lim_{\limnh} \|T_\nh p_\nh u-p_\nh Tu\|_{X_\nh}=0, \qquad& 
  \lim_{\limnh} \|\Bstab_\nh p_\nh u - p_\nh \Bstab u\|_{X_\nh}=0 \qquad \forall u \in X, \\
\label{eq:thm3c}
A_\nh T_\nh = & \Bstab_\nh+\Bcomp_\nh.
\end{align}
\end{subequations}
Then $(A_\nh)_{\nh\in\IN}$ is \regular.
\end{theorem}
\begin{proof}
  Let $(u_\nh)_{\nh\in\IN}$, $u_\nh\in X_\nh$ be a uniformly bounded sequence $\|u_\nh\|_{X_\nh}\leq C$, $(f_\nh)_{\nh\in\IN}$ with $f_\nh := A_\nh u_\nh$ be compact, and $\IN'\subset\IN$ be an arbitrary subsequence.
  Consider a converging subsequence $(f_\nh)_{\nh\in\IN''}$ with $\IN''\subset\IN'$ and denote the limit as $f\in X$ such that $\lim_{n\in\IN''}\|A_\nh u_\nh-p_\nh f\|_{X_\nh}=0$.
  We then obtain from \eqref{eq:thm3c} that $B_\nh T_\nh^{-1} u_\nh +  K_\nh T_\nh^{-1} u_\nh = A_\nh u_\nh = f_n \to f$ for $\nh \in \IN'', \nh \to \infty$.
  Since $T_\nh^{-1}$ is bounded and $(K_\nh)_{\nh\in\IN}$ is \compactop~, $(T_\nh^{-1} u_n)_{\nh \in \IN''}$ is bounded and we can choose a converging subsequence $(g_\nh)_{\nh\in\IN'''}$ with $g_\nh = K_\nh T_\nh^{-1}u_\nh$, $\IN'''\subset\IN''$ and limit $g\in X$ such that $\lim_{\nh\in\IN'''} \| g_\nh -p_\nh g\|_{X_n}=0$.
  We observe that there holds $u_\nh = T_\nh B_\nh^{-1} (f_n - g_n)$.
  Finally, we want to exploit the properties in \eqref{eq:thm3ab} on $T_n$ and $B_n$ to show $\lim_{n\in\IN'''} \|u_\nh-p_\nh T B^{-1}(f-g)\|_{X_h} = 0$ which implies the \compactnessseq of $(u_\nh)_{\nh\in\IN}$.
We start with a triangle inequality
\begin{align*}
    \|u_\nh-&p_\nh T B^{-1}(f-g) \|_{X_\nh}
\leq \underbrace{\|u_\nh-T_\nh B_\nh^{-1}p_\nh(f-g)\|_{X_\nh}}_{I}
+\underbrace{\|p_\nh T B^{-1}(f-g)-T_\nh B_\nh^{-1}p_\nh(f-g)\|_{X_\nh}}_{II}
\end{align*}
and bound the two contributions I and II one after another:
\begin{align*}
 I 
&\leq { \|T_\nh\|_{L(X_\nh)} \|B_\nh^{-1}\|_{L(X_\nh)} \|(f_n-g_n)-p_\nh(f-g)\|_{X_\nh}}
\leq C^2\big( \|f_\nh -p_\nh f\|_{X_\nh} + \|p_\nh-p_\nh g\|_{X_\nh} \big), \\
II &\leq \|p_\nh T B^{-1}(f-g) - T_\nh p_\nh B^{-1}(f-g)\|_{X_\nh}
   + \|T_\nh p_\nh B^{-1}(f-g)-T_\nh B_\nh^{-1}p_\nh(f-g)\|_{X_\nh}\\
&\leq \|p_\nh T B^{-1}(f-g) - T_\nh p_\nh B^{-1}(f-g)\|_{X_\nh}
+ C^2 \|B_\nh p_\nh B^{-1}(f-g)-p_\nh(f-g)\|_{X_\nh},
\end{align*}
where the latter right-hand side terms converge to zero for $n\to\infty$ by the assumptions in \eqref{eq:thm3ab}. Hence $(u_\nh)_{n\in\IN'''}$ converges (to $T B^{-1}(f-g)$) and thus $A_\nh$ is \regular.
\end{proof}

We call a sesquilinear form $a(\cdot,\cdot)$ compact or (weakly) (right $T$-)coercive, if its Riesz representation $A\in L(X)$ (defined by $\inner{Au,u'}_X=a(u,u')$ for all $u,u'\in X$) admits the respective property.

\section{Discrete approximations of the damped time-harmonic Galbrun's equation}\label{sec:galbrun}
In this section we analyze approximations to \eqref{eq:Galbrun}.
After introducing the weak formulation of the problem in \cref{sec:prelim}, we discuss a Helmholtz-type decomposition and a density result in \cref{sec:topo} and \cref{sec:density}, respectively. The discrete approximation is then introduced in \cref{sec:h1conf} and analysed in two steps in \cref{subsec:homo} and \cref{subsec:hetero}, where in \cref{subsec:homo} we treat the case of homogeneous pressure and gravity and treat the general case in \cref{subsec:hetero}.
\subsection{Preliminaries, notation and weak formulation} \label{sec:prelim}
To this end we first set our notation, and specify our assumptions on the parameters and the domain.
Let $\dom\subset\IR^3$ be a bounded Lipschitz polyhedron.
We consider $\dom$ to be the default domain for all functions spaces, i.e.\ $L^2:=L^2(\dom)$, etc..
Let $L^2_0:=\{u\in L^2\colon \mean(u)=0\}$.
Further for a scalar function space $X$ we use the boldface notation for its vectorial variant, i.e.\ $\bX:=(X)^3$.
If not specified otherwise, all function spaces are considered over $\IC$.
We introduce the following subspaces of $\bH^1$ with zero (normal) trace:
\begin{align*}
  \Hz:=\{\u\in \bH^1\colon \nv\cdot\u=0 \text{ on } \partial \dom\} \quad\text{and}\quad
  \Hzz:=(H^1_0)^3,
\end{align*}
where $H^1_0:=\{\u\in H^1\colon \u=0 \text{ on } \partial \dom\}$ is the subspace of $H^1$ with zero trace. 
By $\Cps>0$ we denote the Poincaré-Steklov constant of $\dom$ which satisfies
\begin{align}\label{eq:ConstPS}
  \Cps \|u\|_{H^1} \leq \|\nabla u\|_{\bL^2} \quad\text{for all }u\in H^1_0.
\end{align}

We denote scalar products as $\spl\cdot,\cdot\spr_X$, whereas a scalar product without index always means the $L^2$-scalar product for scalar and vectorial functions.
We employ the notation $A\lesssim B$, if there exists a constant $C>0$ such that $A\leq CB$.
The constant $C>0$ may be different at each occurrence and can depend on the domain $\dom$, the physical parameters $\rho,c_s,p,,\phi,\gamma,\bflow,\omega,\angvel$, and on the sequence of Galerkin spaces $(X_n)_{n\in\IN}$.
However, it will always be independent of the index $n$ and any involved functions which may appear in the terms $A$ and $B$.\\
Let the frequency $\omega\in\IR\mnull$ and the angular velocity of the frame $\angvel\in\IR^3$.
Let the sound speed, density and damping parameter $c_s, \rho, \gamma\colon\dom\to\IR$ be measurable and such that
\begin{align}
\ull{c_s}\leq c_s \leq \ol{c_s}, \qquad
\ull{\rho}\leq \rho \leq \ol{\rho}, \qquad
\ull{\gamma}\leq \gamma \leq \ol{\gamma},
\end{align}
with constants $0<\ull{c_s}, \ol{c_s}, \ull{\rho}, \ol{\rho}, \ull{\gamma}, \ol{\gamma}$.
Let the pressure and gravitational potential $p ,\phi \in W^{2,\infty}$.
Let the source term $\bff\in\bL^2$.
Further let the flow $\bflow\in W^{1,\infty}(\dom,\IR^3)$ such that $\div(\rho\bflow)\in L^2$ and $\nv\cdot\bflow=0$ on $\partial\dom$ and $\bflow$ be compactly supported in $\dom$.
This ensures that the distributional streamline derivative operator $\conv\u:=\bflow\cdot\nabla\u$ is well-defined w.r.t. the inner product $\inner{\rho \cdot,\cdot}$ for $\u\in \bL^2$ \cite{HallaHohage:21}, and we define
\begin{align*}
\X&:=\{\u\in \bL^2\colon \div\u\in L^2,\ \conv \u\in \bL^2, \nv\cdot\u=0 \text{ on } \partial \dom\}
\end{align*}
with inner product
\begin{align*}
    \inner{\u,\u'}_\IX:=\inner{\div \u,\div \u'} + \inner{\conv \u,\conv \u'} + \inner{\u,\u'}
\end{align*}
and the associated norm $\norm{\u}_\IX^2=\inner{\u,\u}_\IX$.
Note that the smoothness $\bflow\in W^{1,\infty}(\dom,\IR^3)$ of the flow will be required to obtain density results for the space $\IX$.
There exists a constant $\Cxh>0$ such that
\begin{align}\label{eq:ConstXH}
\Cxh \|\u\|_\X \leq \|\u\|_{\bH^1} \quad\text{for all }\u\in\bH^1.
\end{align}
We further assume the conservation of mass $\div(\rho\bflow)=0$, which allows us to reformulate \eqref{eq:Galbrun} in the weak form as in \cite{HallaHohage:21}: find $\u\in\Xz$ such that
\begin{equation}
    \acow(\u,\u')=\inner{\bff,\u'}\quad \forall \u'\in\Xz
\end{equation}
with the sesquilinear form
\begin{equations}\label{eq:cowling}
    \acow(\u,\u'):=&\,\inner{c_s^2\rho\div \u,\div \u'} - \inner{\rho\opd \u,\opd \u'}\\
                  &+\inner{\div\u,\grad p\cdot \u'}+\inner{\grad p\cdot\u,\div \u'}+\inner{(\hess(p)-\rho\hess(\phi)) \u,\u'} \\
                  &- i\omega\inner{\gamma\rho\u,\u'}.
\end{equations}

\subsection{Topological decomposition} \label{sec:topo}
A crucial tool to analyse \eqref{eq:cowling} and to construct a proper operator $T$ in \cite{HallaHohage:21} is a Helmholtz-type decomposition of vector fields in $\Xz$.
{To this end let us recall that a vector space $Y$ is called the direct algebraic sum of subspaces $Y_1,\dots,Y_N \subset Y$, denoted by $Y=\bigoplus_{n=1,\dots,N} Y_n$, if each element $y\in Y$ has a unique representation of the form $y=\sum_{n=1}^N y_n$ with $y_n\in Y_n$. We refer to $Y=\bigoplus_{n=1,\dots,N} Y_n$ as the algebraic decomposition of Y.
Note that there exist associated projection operators $P_{Y_n}\colon Y\rightarrow Y_n\colon y\mapsto y_n$ with $\ran P_{Y_n} = Y_n$ and $\ker P_{Y_n} =\bigoplus_{m=1,\dots,N,m\neq n} Y_m$.
An algebraic decomposition of a Hilbert space is called a topological decomposition, if all associated projection operators $P_{Y_n}$ are continuous.
}
We set
\begin{align}
    \begin{split}
&\Vnn:=\{\u\in\Hzz\colon \inner{\nabla\u,\nabla\u'}=0 \text{ for all }\u'\in \Hzz\text{ with }\div\u'=0\},\\
&\Wnn:=\{\u\in\IX\colon \div\u=0\}. \label{eq:Vnn:Wnn}
    \end{split}
\end{align}
Due to \cite[Theorem 4.1]{AcostaDuranMuschietti:06} we know that $\D\v:=\div\v$, $\D\in L(\Vnn,L^2_0)$ is bijective.
We make use of the notation $\D$ for the divergence operator to emphasize that it has a bounded inverse on $\Vnn$, and we will always consider it in the space $\D^{-1}\in L(L^2_0,\Vnn)$.
Note that $\D=\div$ is also bounded and well-defined on $\IX$.
While the choice of $\D$ is deceptively simple in the case of homogeneous pressure and gravity, it is not trivial in the case of heterogeneous pressure and gravity, as we will see in \cref{subsec:hetero}.
On $\Vnn$ the sesquilinear form $\spl\div\cdot,\div\cdot\spr$ defines an inner product equivalent to the $\Hzz$ inner product.

The projections onto $\Vnn$ and $\Wnn$ are given by
\begin{align*}
    {P_V\u:=\D^{-1}\div\u, \qquad P_W\u:=\u-P_V\u.}
\end{align*}
Thus $\Vnn\oplus\Wnn$ is a topological decomposition of $\Xz$.
{If there is no conflict of notation we use the abbreviations $\v:=P_V\u$, $\w:=P_w\u$ for $\u\in\Xz$.}

\subsection{Density results}\label{sec:density}

{
\begin{proposition}[Variation of Prop.\ 3.5 of \cite{Bredies:08}]
\label{prop:Bredies}
Let $l\in\mathbb{N}$ and $\Lambda\colon \mathcal{D}(\Lambda)\subset L^2(\dom,\mathbb{R}^3)\to L^2(\dom,\mathbb{R}^l)$ with $C_0^\infty(\dom,\mathbb{R}^3)\subset\mathcal{D}(\Lambda)$ be a closed linear operator with the property that
\begin{enumerate}
	\item $u\in\mathcal{D}(\Lambda)$ if and only if for each $\zeta\in C_0^\infty(\dom)$ follows $\zeta u\in\mathcal{D}(\Lambda)$,
	\item for each $u\in\mathcal{D}(\Lambda)$ and $\zeta\in C_0^\infty(\dom)$ follows $\supp \Lambda(\zeta u)\subset \supp \zeta$,
	\item for each $u\in\mathcal{D}(\Lambda)$ with compact support in $\dom$, there exists a $\delta_0>0$ such that the sequence of mollified $u_\delta:=u*G_\delta$ satisfies $\|\Lambda u_\delta\|_{L^2(\dom,\mathbb{R}^l)} \leq C$ for every $\delta\in(0,\delta_0)$ and some $C>0$.
\end{enumerate}
Then, for each $\epsilon>0$ and $u\in L^2(\dom,\mathbb{R}^3)$ with $\Lambda u\in L^2(\dom,\mathbb{R}^l)$, there exists a $\tilde u\in C^\infty(\dom,\mathbb{R}^3)$ such that
\begin{align}
	\|u-\tilde u\|_{L^2(\dom,\mathbb{R}^3)}^2
	+\|\Lambda u-\Lambda \tilde u\|_{L^2(\dom,\mathbb{R}^l)}^2
	<\epsilon.
	\label{eq:Bredies}
\end{align}
\end{proposition}
}

\begin{theorem}\label{thm:densXH}
Let $\bflow\in W^{1,\infty}(\dom,\IR^3)$ and $\bflow$ be compactly supported in $\dom$.
Then $\Hzz$ is dense in $\Xz$.
\end{theorem}
\begin{proof}

Let $\chi \in C_0^\infty(\dom)$ be a cut-off function with values in $[0,1]$, $\dom \supset G_2 \supset G_1 \supset\supp \bflow$,
$\chi=1$ on $G_1$ and $\chi=0$ on $\dom\setminus G_2$,
where $\operatorname{dist}(\partial G_2, \partial \dom) > 0$.
Let $\u\in\IX$ and $\epsilon>0$.
Since $\|(1-\chi)\u\|_{\IX}=\|(1-\chi)\u\|_{H(\div;\dom)}$ we can find $\tilde\u_1\in\mathbf{C}_0^\infty(\dom)$ such that $\|(1-\chi)\u-\tilde \u_1\|_{\IX}<\epsilon/2$, see, e.g., \cite{ErnGuermond_FE_I}.
To find a suitable smooth approximation of $\chi\u$ we apply \Cref{prop:Bredies} to $\Lambda\u:=(\div\u,\conv\u)^\top$, i.e., $l=4$.
The first assumption of \Cref{prop:Bredies} follows from the product rule (see, e.g., \cite[Lem.~3.7]{Bredies:08} for details on $\conv$).
The second assumption of \Cref{prop:Bredies} holds, because $\Lambda$ is a differential operator.
The third assumption of \Cref{prop:Bredies} follows from \cite[Lem.~3.8]{Bredies:08} and convenient manipulations for the smoothing in $H(\div)$, i.e.,
$\div(\u*G_\delta)(\bx)=\int_\dom \u(\by) \cdot \nabla_\bx G_\delta(\bx-\by)d\by
=-\int_\dom \u(\by) \cdot \nabla_\by G_\delta(\bx-\by)d\by
=\int_\dom G_\delta(\bx-\by) \div_\by \u(\by) d\by$.
The claimed bound follows now from the properties of $G_\delta$.
Thus there exists $\tilde\u_2\in\mathbf{C}^\infty(\dom)$ such that $\|\chi\u-\tilde \u_2\|_{\IX}<\epsilon/2$.
Since the support of $\chi\u$ is compact in $\dom$, $\tilde\u_2$ can be choosen with compact support too and hence satisfies the necessary boundary condition.
Thus the proof is finished.
\end{proof}

\begin{theorem}\label{thm:densXC}
Let $\bflow\in W^{1,\infty}(\dom,\IR^3)$ and $\supp\bflow$ be compact in $\dom$.
Then $\bC_0^\infty$ is dense in $\Xz$.
\end{theorem}
\begin{proof}
Since $\bC^\infty_0$ is dense in $\Hzz$ the claim follows from Theorem~\ref{thm:densXH} and \eqref{eq:ConstXH}.
\end{proof}

\subsection{\texorpdfstring{$\bH^1$}{H1}-conforming discretization}\label{sec:h1conf}

Let $(\calT_n)_{n\in\IN}$ be a sequence of shape-regular simplical meshes of $\dom$ with maximal element diameter $h_n\to0$ for $n\to\infty$.
For $k\in\IN$ we denote by $P_k$ the space of scalar polynomials with maximal degree $k$.
We consider finite element spaces
\begin{align*}
  \Xn:=\{\u\in\Hz \colon \u|_T \in (P_k(T))^3 ~ \forall T\in\calT_n\},
\end{align*}
with fixed uniform polynomial degree $k\in\IN$.
It readily follows $\Xnz\subset\Hz\subset\Xz$.

Let us note that the previous assumption that $\dom$ is polygonal is crucial for $\Xn$ to be a proper finite element space with the usual approximation quality. We discuss the construction of such a finite element space in \Cref{sec:Xn:construction}. For curved boundaries, especially in the case of curved boundaries that are approximated with only $C^0$-continuous discrete boundaries, the construction of $\Xn$ is hardly possible or computationally unfeasible. In these cases, one typically resorts to Lagrange multiplier-based or Nitsche-like techniques in order to weakly impose the boundary condition $\bu \cdot \nv=0$ through the variational formulation that is then posed on 
$
\{\u\in\bH^1 \colon \u|_T \in (P_k(T))^3 ~ \forall T\in\calT_n\}.
$
In the numerical examples below we will use a Nitsche-based (weak) imposition of the boundary conditions while in the analysis we assume $\bu \cdot \nv=0$ to be imposed as \emph{essential} boundary conditions in $\Xn$.

$\Xn$ allows for proper approximation of $\bu \in \Xz$:
\begin{lemma}\label{lem:XnDense}
It holds
\begin{align*}
\lim_{n\to\infty} \inf_{\u_n'\in\Xnz} \|\u-\u_n'\|_{\IX} =0 \quad\text{for each}\quad \u\in\Xz.
\end{align*}
\end{lemma}
\begin{proof}
Let $\bu\in\Xz$ be given.
Since $\bC^\infty_0$ is dense in $\Xz$ (see Thm.~\ref{thm:densXC}) we can find for each $\epsilon>0$ a function $\u_\epsilon\in\bC^\infty_0$ such that $\|\u-\u_\epsilon\|_\IX<\epsilon$.
Further, the canonical interpolation operator $I_{h_n}$ is well defined for $\u_\epsilon$ and yields the estimate $\|\u_\epsilon- I_{h_n}\u_\epsilon\|_{\bH^1} \leq C h_n \|\u_\epsilon\|_{\bH^2}$ with a constant $C>0$ independent of ${h_n}$.
Since $\u_\epsilon$ has compact support it also follows that $I_{h_n}\u_\epsilon\in\Hzz$ and thus $I_{h_n}\u_\epsilon\in\Xnz$.
Hence we estimate
\begin{align*}
  \lim_{n\to\infty} \inf_{\u_n'\in\Xnz} \|\u-\u_n'\|_{\IX}
  &\leq \lim_{n\to\infty} \inf_{\u_n'\in\Xnz} (\|\u-\u_\epsilon\|_\IX+\|\u_\epsilon-\u_n'\|_{\IX})
  \leq \epsilon + \lim_{n\to\infty} \|\u_\epsilon-I_{h_n}\u_\epsilon\|_{\IX} \\
  &\lesssim \epsilon + \lim_{n\to\infty} \|\u_\epsilon-I_{h_n}\u_\epsilon\|_{\bH^1} 
  \lesssim \epsilon + \lim_{n\to\infty} h_n\|\u_\epsilon\|_{\bH^2} 
= \epsilon.
\end{align*}
Since $\epsilon>0$ was chosen arbitrarily it follows $\lim_{n\to\infty} \inf_{\u_n'\in\Xnz} \|\u-\u_n'\|_{\IX}=0$.
\end{proof}
Let $\PXnz \in L(\Xz,\Xnz)$ be the $\IX$-orthogonal projection onto $\Xnz$.
Lemma~\ref{lem:XnDense} implies that $\lim_{n\to\infty} \|\u-\PXnz\u\|_\X=0$ for each $\u\in\Xz$.

Based on $\Xn$ we can formulate the discrete problem as:
\begin{align}\label{eq:H1disc}
    \begin{split}
        &\text{find } \u_n\in \Xnz\text{ s.t.}%
        \quad \acow(\u_n,\u'_n)
        =\inner{\bff,\u'_n}\quad \forall \u'_n\in \Xnz.
    \end{split}
\end{align}
Let $A\in L(\Xz)$ be the operator associated to $\acow(\cdot,\cdot)$ and $A_n:=\PXnz A|_{\Xnz} \in L(\Xnz)$.
Then the introduced Galerkin approximation constitutes a discrete approximation scheme as described in Section~\ref{subsec:DAS}, whereat $p_n=\PXnz$.
To guarantee the stability of the approximations we impose the following assumption.
Let
\begin{align} \label{eq:Qn}
\Qn:=\{f\in L^2_0\colon f|_T\in P_{k-1}(T) \quad\forall T\in\mesh\}
\end{align}
and let $P_{Q_n}\in L(L^2_0,Q_n)$ be the associated orthogonal projection.

A key observation of the following analysis is that a discrete inf-sup-stability for the discrete divergence operator and the spaces $\Xn$ and $Q_n$ allows to obtain a discrete counterpart of the Helmholtz-type decomposition that is required for the discrete $T_n$ operator in the $T$-coercivity analysis.
\begin{assumption}\label{ass:DivStable}
    There exists a constant ${\betad}>0$ such that
$$\inf_{f_n\in\Qn\setminus\{0\}} \sup_{\u_n\in\Xn\setminus\{0\}} \frac{|\inner{\div\u_n,f_n}|}{\|\nabla\u_n\|_{(L^2)^{3\times3}} \|f_n\|_{L^2}}>\betad$$ for all $n\in\IN$.
\end{assumption}
The choice of $\Qn$ in \eqref{eq:Qn} relates to Scott-Vogelius elements in the discretization of the Stokes problem.
In order to ensure its stability and hence to make sure that \cref{ass:DivStable} is satisfied it is usually necessary to apply special meshes (barycentric refinement) and/or sufficiently large polynomial degree $k$, see e.g.\ \cite{Neilan:20} and \cite{ScottVogelius:85,Zhang:11a,ArnoldQin:92,Zhang:05,GuzmanNeilan:18}.

While the Scott–Vogelius element satisfies the stronger condition $\div(\IX_n) \subset \Qn$, this property is not essential for the validity of the analysis, as will also be clarified in \cref{rem:th}.
\cref{ass:DivStable} can often be relaxed if $\Qn$ is replaced by another finite element space and in the discrete formulation $\div$ is replaced by $\div_h := P_{\Qn} \div$.
We will comment on this type of discretizations and the necessary adjustments in the analysis in more detail in \cref{rem:th} after the first a priori error bounds, below.

\subsection{Homogeneous pressure and gravity}\label{subsec:homo}

In this section we consider a simplified case of \eqref{eq:Galbrun} in which the pressure and gravitational potential are assumed to be constant before we consider the general case in the subsequent section.
\eqref{eq:cowling} reduces to
\begin{align*}
    \acow(\u,\u')\!=\!\inner{c_s^2\rho\div \u,\div \u'}\! -\! \inner{\rho\opd \u,\opd \u'}\! -\! i\omega\inner{\gamma\rho\u,\u'}.
\end{align*}
We aim to establish the stability of $(A_n)_{n\in\IN}$ by means of \Cref{thm:Tcomp} and \Cref{lem:stable}.
To this end we need to construct operators $\Tn$ with respective properties.
Of course the natural approach is to mimic the analysis from the continuous level \cite{HallaHohage:21}.
However, for the analysis in this article we will rely on a slightly different construction than used in \cite{HallaHohage:21}.
The reason thereof is that this new variant can be mimicked more easily on the discrete level.
While the analysis presented here is an important setup for the discrete problem, compared to the results in \cite{HallaHohage:21} it is suboptimal, as the assumption on the Mach number $\|c_s^{-1}\bflow\|_{\bL^\infty}$ is more restrictive.
\begin{lemma}\label{lem:wTc}
Let $\beta>0$ be the $\inf$-$\sup$ constant of the divergence on $\dom$.
Let $\|c_s^{-1}\bflow\|_{\bL^\infty}^2<\beta^2 \frac{\ull{c_s}^2\ull{\rho}}{\ol{c_s}^2\ol{\rho}}$.
Let $T:=\PV-\PW$. Then {$T\in L(\Xz)$ is bijective with inverse $T^{-1}=T$ and} $A$ is weakly right $T$-coercive.
\end{lemma}
\begin{proof}
{Since $P_V, P_W$ are the projections of a topological decomposition it holds that $T\in L(\Xz)$ and $TT=(P_V-P_W)(P_V-P_W)=P_VP_V+P_WP_W=P_V+P_W=I$}.
Using $T\u=\v-\w$ we have that $\spl AT\u,\u \spr_\IX=\acow(\v-\w,\v+\w)$.
It then holds $AT=B+K$ for $B$ and $K$ defined by
\begin{align}
\spl B\u,\u' \spr_\IX:= &
\inner{c_s^2\rho\div \v,\div \v'}
-\inner{\rho i\partial_\bflow\v,i\partial_\bflow \v'} \nonumber \\
&-\inner{\rho i\partial_\bflow\v,\opd \w'}+\inner{\rho \opd \w,i\partial_\bflow\v'} \label{eq:B} \\
&+\inner{\rho\opd \w,\opd \w'}
+i\omega\inner{\gamma\rho\w,\w'} \nonumber \\
\spl K\u,\u' \spr_\IX:= & \nonumber 
-\inner{\rho (\omega+i\angvel\times) \v, (\omega+i\angvel\times) \v'}
-\inner{\rho (\omega+i\angvel\times) \v, i\partial_\bflow \v'}
-\inner{\rho i\partial_\bflow \v, (\omega+i\angvel\times) \v'}\\ \label{eq:K}
&-\inner{\rho (\omega+i\angvel\times) \v, \opd \w'}
+\inner{\rho \opd \w, (\omega+i\angvel\times) \v'}\\ \nonumber
&-i\omega \spl \rho\gamma\v,\v'\spr
+i\omega \spl \rho\gamma\w,\v'\spr
-i\omega \spl \rho\gamma\v,\w'\spr, ~~\text{
  for all}~ \u,\u'\in\IX.
\end{align}
{The terms appearing in definition of $K$ can be represented e.g.\ as
\begin{align*}
\inner{\rho (\omega+i\angvel\times) \v, \opd \w'}=\inner{P_W^* B_{\opd}^* M_{\rho (\omega+i\angvel\times)}E_{V,\bL^2}P_V\u,\u'}_{\Xz}
\end{align*}
with the embedding $E_{V,\bL^2}\in L(V,\bL^2)$, the multiplication operator $M_{\rho (\omega+i\angvel\times)}\in L(\bL^2)$ and $B_{\opd}\in L(W,\bL^2)$, $B_{\opd}\w:=\opd\w$.
Since the embedding $\bH^1\hookrightarrow\bL^2$ is compact, $V$ embeds continously into $\bH^1$ and each term in \eqref{eq:K} contains at least one operator $E_{V,\bL^2}$ or $E_{V,\bL^2}^*$ it follows that $K\in L(\Xz)$ is compact.
}
We now show that $B$ is coercive and hence that $A$ is bijective. Let $\tau\in(0,\pi/2)$.
We compute
\begin{align*}
\frac{1}{\cos\tau}\Re \left( e^{-i\tau\sign\omega} \spl B\u,\u \spr_\IX\right)
=&\inner{c_s^2\rho\div \v,\div \v}
-\inner{\rho i\partial_\bflow\v,i\partial_\bflow \v}
+\tan\tau |\omega|\inner{\gamma\rho\w,\w}
\\
&+\inner{\rho\opd \w,\opd \w}\\
&-2\tan\tau\sign\omega \Im\left( \inner{\rho i\partial_\bflow\v,\opd \w} \right)
\end{align*}
We estimate the last term by the Cauchy--Schwarz inequality and the weighted Young inequality
$|2ab| \leq (1-\epsilon)^{-1}a^2 + (1-\epsilon)b^2$ with an additional parameter $\epsilon\in(0,1)$,
$a = \tan \tau  \| \sqrt{\rho} \partial_\bflow\v\|_{L^2}$, and $b=\|\sqrt{\rho}\opd\w\|_{L^2}$ and obtain
\begin{align*}
\frac{1}{\cos\tau}\Re \left( e^{-i\tau\sign\omega} \spl B\u,\u \spr_\IX\right)
&\geq\inner{c_s^2\rho\div \v,\div \v}
-\big(1+(1-\epsilon)^{-1}\tan^2\tau\big)\inner{\rho i\partial_\bflow\v,i\partial_\bflow \v}\\
&\qquad+\epsilon\inner{\rho\opd \w,\opd \w}+\tan\tau |\omega|\inner{\gamma\rho\w,\w}
\end{align*}
We estimate further
\begin{align*}
\inner{c_s^2\rho\div \v,\div \v}
&-\big(1+(1-\epsilon)^{-1}\tan^2\tau\big)\inner{\rho i\partial_\bflow\v,i\partial_\bflow \v}\\
&\geq \left(\beta^2 \ull{c_s}^2 \ull{\rho} - \ol{c_s}^2\ol{\rho} \|c_s^{-1}\bflow\|_{\bL^\infty}^2 \big(1+(1-\epsilon)^{-1}\tan^2\tau\big) \right) \inner{\nabla \v,\nabla \v}.
\end{align*}
Due to the assumption of this lemma we can choose small enough $\tau\in(0,\pi/2)$ and $\epsilon>0$ such that the constant in the right hand-side is positive.
Since $\|\v\|_\IX \lesssim \|\nabla\v\|_{L^2}$ for $\v\in\Hzz$ this yields coercivity in $\v$.
As in \cite{HallaHohage:21} a weighted Young's inequality shows that
\begin{align*}
\epsilon\inner{\rho\opd \w,\opd \w}+\tan\tau |\omega|&\inner{\gamma\rho\w,\w}
\gtrsim \norm{\partial_\bflow\w}_{L^2}^2+\norm{\w}_{L^2}^2
=\norm{\w}_\IX^2.
\end{align*}
Thus $|\spl B\u,\u \spr_\IX|\gtrsim \norm{\v}_\IX^2+\norm{\w}_\IX^2\gtrsim \norm{\u}_\IX^2$ and the claim follows.
\end{proof}

\subsubsection{Regular approximation}\label{subsubsec:RegApr}

{
\begin{lemma}
\label{lem:VnWn}
Let Assumption~\ref{ass:DivStable} be satisfied.
Then the spaces
\begin{align} \begin{split}\label{eq:Vhnn:Whnn}
&\Vhnn:=\{\u_n\in \Xn\cap\Hzz\colon \inner{\nabla\u_n,\nabla\u'_n}=0\quad \forall \u'_n\in \Hzz\cap\Whnn\}, \\
&\Whnn:=\{\u_n\in \Xn\colon \div \u_n=0\},
\end{split}\end{align}
form a topological decomposition of $\Xn$ with projections
\begin{align*}
	P_{\Vhnn}\u_n:=\D_n^{-1}\div\u_n, \qquad P_{\Whnn}\u_n:=\u_n-P_{\Vhnn}\u_n,
\end{align*}
being uniformly bounded in $n\in\mathbb{N}$, where for $q_n\in Q_n$ the function $\D_n^{-1}q_n\in\Vhnn$ is the unique solution to
\begin{align}\label{eq:vn}
	\text{find }\v_n\in\Vhnn\text{ such that } \div \v_n = q_n, %
\end{align}
i.e.\ $\D_n^{-1}\in L(Q_n,\Vhnn)$.
\end{lemma}
\begin{proof}
Assumption~\ref{ass:DivStable} ensures that \eqref{eq:vn} admits a unique solution $\v_n$ which satisfies $\betad\|\nabla\v_n\|_{(L^2)^{3\times3}}\leq\|q_n\|_{L^2}$.
Since $P_{\Vhnn}P_{\Vhnn}\u_n=\D_n^{-1}\div P_{\Vhnn}\u_n=\D_n^{-1}\div \u_n=P_{\Vhnn}\u_n$, $P_{\Vhnn}$ is indeed a projection.
Assumption~\ref{ass:DivStable} ensures that $P_{\Vhnn}$ is uniformly bounded.
Due $\ker P_{\Vhnn}=\Whnn$ the spaces $\Vhnn,\Whnn$ form indeed a topological decomposition of $\Xn$.
\end{proof}
We abbreviate $\v_n:=P_{V_n}\u_n$, $\w_n:=P_{W_n}\u_n$ for $\u_n\in\Xn$.
\begin{lemma}
\label{lem:Tn}
Let Assumption~\ref{ass:DivStable} be satisfied.
Then $$T_n:=P_{\Vhnn}-P_{\Whnn}=T_n^{-1}\in L(\Xn)$$ is uniformly bounded in $n\in\mathbb{N}$.
\end{lemma}
\begin{proof}
Since the spaces $\Vhnn,\Whnn$ form a topological decomposition of $\Xn$ it follows that $T_n=T_n^{-1}$.
The uniform boundedness of $T_n,T_n^{-1}$ follow from the uniform boundedness of $P_{\Vhnn}$.
\end{proof}
}
\begin{lemma}\label{lem:DnInvPQnDiv}
For each $\v\in\Vnn$ it holds that $\lim_{n\to\infty}\norm{\v-\D_n^{-1}\PQn\div\v
}_{\bH^1}=0$.
\end{lemma}
\begin{proof}
$(\v,0)$ solves the problem to find $(\u,\z)\in \Hzz\times (\Wnn\cap\Hzz)$ such that
\begin{align*}
\inner{\div \u,\div \u'}
+\inner{\nabla\u,\nabla\z'}
+\inner{\nabla\z,\nabla\u'}
= \inner{\div \v,\div \u'},
\end{align*}
for all $(\u',\z')\in \Hzz\times(\Wnn\cap\Hzz),$
and $(\D_n^{-1}\PQn\div\v,0)$ solves the problem to find $(\u_n,\z_n)\in \Xn\cap\Hzz\times (\Whnn\cap\Hzz)$ such that
\begin{align*}
\inner{\div \u_n,\div \u'_n}
+\inner{\nabla\u_n,\nabla\z'_n}
+\inner{\nabla\z_n,\nabla\u'_n}
= \inner{\div \v,\div \u'_n},
\end{align*}
for all $(\u_n',\z_n')\in \Xn\cap\Hzz\times(\Whnn\cap\Hzz)$.
The latter is a conforming Galerkin approximation of the former.
It can be seen that both equations are uniformly stable by testing with $(\v+\z,\w)$ and $(\v_n+\z_n,\w_n)$ respectively.
With a Céa lemma, it only remains to show that $\lim_{n\to\infty}\inf_{\u_n\in\Xn\cap\Hzz}\norm{\u-\u_n}_{\bH^1}=0$ and $\lim_{n\to\infty}\inf_{\z_n\in\Whnn\cap\Hzz}\norm{\z-\z_n}_{\bH^1}=0$ for $\u=\v\in \Vnn \subset \Hzz$ and $\z=0$ respectively.
The first result is standard while the second is trivial as $\z = 0 \in \Whnn$.
\end{proof}
Next we shall establish the point-wise limit of $T_n$.
\begin{lemma}\label{lem:TnLimit}
For each $\u\in\IX$ it holds $\lim_{n\to\infty}\|T_n\PXnz\u-\PXnz T\u\|_\IX=0$.
\end{lemma}
\begin{proof}
It suffices to prove $\lim_{n\to\infty}\|P_{\Vhnn}\PXnz\u-\PXnz P_{\Vnn}\u\|_\IX=0$.
We use $P_{\Vhnn}=\D_n^{-1}\PQn\div$ and estimate
\begin{align*}
&\|P_{\Vhnn}\PXnz\u-\PXnz P_{\Vnn}\u\|_\IX
\leq \|\D_n^{-1}\PQn\div\PXnz\u-P_{\Vnn}\u\|_\IX+\|P_{\Vnn}\u-\PXnz P_{\Vnn}\u\|_\IX \\
&\!\leq \|\D_n^{-1}\PQn\div\u\!-\!P_{\Vnn}\u\|_\IX\|+\|\D_n^{-1}\PQn\div\|_{L(\IX)}\|\u\!-\!\PXnz\u\|_\IX+\|P_{\Vnn}\u\!-\!\PXnz P_{\Vnn}\u\|_\IX.
\end{align*}
The claim follows now from $\div\u=\div P_{\Vnn}\u$, Lemma~\ref{lem:DnInvPQnDiv},
\eqref{eq:ConstXH} and the point-wise convergence of $\PXnz$ (see Lemma~\ref{lem:XnDense}).
\end{proof}

\begin{lemma}\label{lem:reg}
If $\|c_s^{-1}\bflow\|_{\bL^\infty}^2<\betad^2 \frac{\ull{c_s}^2\ull{\rho}}{\ol{c_s}^2\ol{\rho}}$, then $(A_n)_{n\in\IN}$ is \regular,\ in the sense of \cref{def:DAS}.
\end{lemma}
\begin{proof}
We apply \Cref{thm:Tcomp}.
In the previous part of this Section~\ref{subsubsec:RegApr} we already constructed $T_n$ and showed that $T_n\in L(\Xn)$ and $T_n^{-1}=T_n\in L(\Xn)$ are uniformly bounded.
Further, \Cref{lem:TnLimit} shows that $T_n$ converges pointwise.
Next we need to split $A_nT_n=B_n+K_n$ into a stable part $B_n\in L(\Xn)$ and a compact part $K_n\in L(\Xn)$.
To do so we stick very closely to the lines of \cite{HallaHohage:21}.
Recall that $\spl A_nT_n\u_n,\u_n \spr_\IX=\acow(\v_n-\w_n,\v_n+\w_n)$.
Hence it holds $A_nT_n=B_n+K_n$ with $B_n$ and $K_n$ defined by
\begin{align*}
\spl B_n\u_n,\u'_n \spr_\IX:=&
\inner{c_s^2\rho\div \v_n,\div \v_n'}
-\inner{\rho i\partial_\bflow\v_n,i\partial_\bflow \v_n'}\\
&-\inner{\rho i\partial_\bflow\v_n,\opd \w_n'}+\inner{\rho \opd \w_n,i\partial_\bflow\v_n'}\\
&+\inner{\rho\opd \w_n,\opd \w_n'}
+i\omega\inner{\gamma\rho\w_n,\w_n'}
\end{align*}
and
\begin{align*}
    \spl K_n\u_n,\u'_n \spr_\IX:=
&-\inner{\rho (\omega+i\angvel\times) \v_n, (\omega+i\angvel\times) \v_n'}
 -\inner{\rho (\omega+i\angvel\times) \v_n, i\partial_\bflow \v_n'}
-\inner{\rho i\partial_\bflow \v_n, (\omega+i\angvel\times) \v_n'}\\
&-\inner{\rho (\omega+i\angvel\times) \v_n, \opd \w_n'}
+\inner{\rho \opd \w_n, (\omega+i\angvel\times) \v_n'}\\
&-i\omega \spl \rho\gamma\v_n,\v_n'\spr
+i\omega \spl \rho\gamma\w_n,\v_n'\spr
-i\omega \spl \rho\gamma\v_n,\w_n'\spr
\end{align*}
for all $\u_n,\u_n'\in\Xn$.
The operator $K_n$ is compact due to the compact Sobolev embedding from $\Vhnn\subset\bH^1$ to $\bL^2$.
It is straightforward to see that $B_n$ is uniformly bounded
and that $B_n$ converges pointwise to the operator $B\in L(\IX)$ defined in 
\eqref{eq:B}.
The uniform coercivity of $B_n$ follows along the lines of the proof of \Cref{lem:wTc}, with the constant $\beta$ replaced by $\betad$.
Hence the claim is proven.
\begin{extproof}
Let $\tau\in(0,\pi/2)$.
We compute
\begin{align*}
\frac{1}{\cos\tau}\Re \left( e^{-i\tau\sign\omega} \spl B_n\u_n,\u_n \spr_\IX\right)
&=\inner{c_s^2\rho\div \v_n,\div \v_n}
-\inner{\rho i\partial_\bflow\v_n,i\partial_\bflow \v_n}\\
&+\inner{\rho\opd \w_n,\opd \w_n}+\tan\tau |\omega|\inner{\gamma\rho\w_n,\w_n}\\
&-2\tan\tau\sign\omega \Im\left( \inner{\rho i\partial_\bflow\v_n,\opd \w_n} \right)
\end{align*}
We estimate the last term by the Cauchy--Schwarz inequality and the weighted Young inequality
$|2ab| \leq (1-\epsilon)^{-1}a^2 + (1-\epsilon)b^2$ with an additional parameter $\epsilon\in(0,1)$,
$a = \tan \tau  \| \sqrt{\rho} \partial_\bflow\v_n\|_{L^2}$, and $b=\|\sqrt{\rho}\opd\w_n\|_{L^2}$ and obtain
\begin{align*}
\frac{1}{\cos\tau}\Re \left( e^{-i\tau\sign\omega} \spl B_n\u_n,\u_n \spr_\IX\right)
&\geq\inner{c_s^2\rho\div \v_n,\div \v_n}
-\big(1+(1-\epsilon)^{-1}\tan^2\tau\big)\inner{\rho i\partial_\bflow\v_n,i\partial_\bflow \v_n}\\
&+\epsilon\inner{\rho\opd \w_n,\opd \w_n}+\tan\tau |\omega|\inner{\gamma\rho\w_n,\w_n}
\end{align*}
We estimate further
\begin{align*}
\inner{c_s^2\rho\div \v_n,\div \v_n}
&-\big(1+(1-\epsilon)^{-1}\tan^2\tau\big)\inner{\rho i\partial_\bflow\v_n,i\partial_\bflow \v_n}\\
&\geq \left(\betad \ull{c_s}^2 \ull{\rho} - \ol{c_s}^2\ol{\rho} \|c_s^{-1}\bflow\|_{\bL^\infty}^2 \big(1+(1-\epsilon)^{-1}\tan^2\tau\big) \right) \inner{\nabla \v_n,\nabla \v_n}.
\end{align*}
Due to the assumption of this lemma we can choose small enough $\tau\in(0,\pi/2)$ and $\epsilon>0$ such that the constant in the right hand-side is positive.
Since $\|\v\|_\IX \lesssim \|\nabla\v\|_{L^2}$ for $\v\in\Hzz$ this yields coercivity in $\v_n$.
As in \cite{HallaHohage:21} a weighted Young's inequality shows that
\begin{align*}
\epsilon\inner{\rho\opd \w,\opd \w}+\tan\tau |\omega|\inner{\gamma\rho\w,\w}
\gtrsim \norm{\partial_\bflow\w_n}_{L^2}^2+\norm{\w_n}_{L^2}^2
=\norm{\w_n}_\IX^2.
\end{align*}
Thus $|\spl B_n\u_n,\u_n \spr_\IX|\gtrsim \norm{\v_n}_\IX^2+\norm{\w_n}_\IX^2\gtrsim \norm{\u_n}_\IX^2$ and the claim follows.
\end{extproof}
\end{proof}

\subsubsection{Convergence}

\begin{theorem}\label{thm:conv}
Let $p$ and $\phi$ be constant.
Let $\u$ be the solution to~\eqref{eq:Galbrun}.
Let Assumption~\ref{ass:DivStable} be satisfied and $\|c_s^{-1}\bflow\|_{\bL^\infty}^2<\betad^2 \frac{\ull{c_s}^2\ull{\rho}}{\ol{c_s}^2\ol{\rho}}$.
Then there exists an index $n_0>0$ such that for all $n>n_0$ the solution $\u_n$ to \eqref{eq:H1disc} exists
and $\u_n$ converges to $\u$ in the $\IX$-norm with the best approximation estimate
$\|\u-\u_n\|_\IX \lesssim \inf_{\u_n'\in\Xn} \|\u-\u_n'\|_\IX$.
\end{theorem}
\begin{proof}
Due to Lemma~\ref{lem:reg} the approximation scheme $(A_n)_{n\in\IN}$ is regular.
Since $A$ is bijective~\cite{HallaHohage:21} the claim follows from \Cref{lem:stable,lem:conv}.
\end{proof}
\begin{remark}\label{rem:convrate}
Note that for smooth solutions $\bu\in \bH^{1+s}$, $s>0$ we can obtain convergence rates by convenient techniques:
\begin{align*}
\inf_{\u_n'\in\Xn} \|\u-\u_n'\|_\IX
\lesssim \inf_{\u_n'\in\Xn} \|\u-\u_n'\|_{\bH^1}
\lesssim h^{\min(s,k)} \|\bu\|_{\bH^{1+s}}.
\end{align*}
\end{remark}

\begin{remark}\label{rem:th}
    The considered discrete setting can be generalized by replacing the divergence operator $\operatorname{div}$ in the discrete formulation by a discrete version $\operatorname{div}_h: \Xn \to Q_n$ with a space $Q_n$ that is potentially different to the one in \eqref{eq:Qn}. In this case also the \cref{ass:DivStable} would be relaxed w.r.t.  $\operatorname{div}$ and $Q_n$. One important case which is known as the \emph{Taylor-Hood} discretization in fluid dynamics is obtained from $Q_n = \{f\in L^2_0\colon f|_T\in P_{k-1}(T) \quad\forall T\in\mesh\} \cap H^1$ and $\operatorname{div}_h = P_{Q_n} \operatorname{div}$. For the implementation of $\operatorname{div}_h$ one typically introduces an auxiliary variable, the so-called pseudo-pressure so that
    $\inner{\operatorname{div}_h \bu_n, \operatorname{div}_h \bu_n'}$ becomes
    $\inner{q_n,\div \bu_n'}+\inner{\operatorname{div} \bu_n, q_n'}-\inner{q_n, q_n'}$ where $\bu_n$ and $q_n$ and $\bu_n'$ and $q_n'$ are the trial and the test functions in $\Xn$ and $Q_n$, respectively.
    Let us briefly sketch the changes in the analysis that would be necessary to account for this change in the discrete formulation. First, note that replacing $\operatorname{div}$ with $\div_h$ in \eqref{eq:H1disc} would lead to a \emph{non-conforming} discretization. Hence, we would need to prove \asymptoticconsistency, i.e. that the corresponding sequence of discrete operators $A_n$ \approximates~ $A$ which has been trivial for the Galerkin approximation. In the discrete subspace splitting  $\Whnn$ would need to be defined w.r.t. to $\div_h$ (instead of $\div$) as well as the corresponding projection onto $\Vhnn$ in \eqref{eq:vn}. With only minor changes also the proof of \cref{lem:DnInvPQnDiv} would carry over to this setting so that finally convergence of the corresponding discrete solution $\bu_n$ to the continuous solution $\bu$ would follow. Alternatively, an equivalent conforming discretization could be analysed by introducing the pseudo-pressure formulation already on the continuous level. 
    In the remainder of the analysis in this manuscript we will continue to focus on to the case of the divergence operator $\operatorname{div}$ and the space $Q_n$ as in \eqref{eq:Qn}. However, in the numerical examples below we will also consider a Taylor-Hood-type discretization and compare it with the chosen setting of Scott-Vogelius-type elements.
  \end{remark}

\subsection{Heterogeneous pressure and gravity}\label{subsec:hetero}

In this section we expand the analysis from the previous section and consider heterogeneous pressure $p$ and gravitational potential $\phi$.

\subsubsection{Analysis on the continuous level}

As in \cite{HallaHohage:21} we introduce $\q:=c_s^{-2}\rho^{-1}\nabla p$ and express
\begin{align}\label{eq:aq}
\acow(\u,\u')&=\inner{c_s^2\rho(\div + \q\cdot)\u,(\div + \q\cdot)\u'}
-\inner{\rho\opd \u,\opd\u'}\nonumber \\
&-i\omega \inner{\rho\gamma\u,\u'}
+\inner{(\hess(p)-\rho\hess(\phi)-c_s^2\rho \, \q\otimes \q)\u,\u'}.
\end{align}
However, in the forthcoming analysis we will deviate from \cite{HallaHohage:21} and avoid the introduction of an additional third space $Z$ in the topological decomposition of $\IX$.
Consider now the divergence operator $\D\in L(\Vnn,L^2_0), \D \v:=\div \v$.
We know that $\D^{-1}\in L(L^2_0,\Vnn)$.
For heterogeneous pressure our analysis leads us to consider $\D\bv+\q\cdot\bv$ instead of $\D\bv$.
A necessary ingredient for our analysis is that the new operator $\D+\q\cdot$ is invertible on suitable spaces.
Since we cannot ensure this property for $\D+\q\cdot$, we work instead with a slight modification.
\begin{lemma}
\label{lem:Dt}
There exist operators $\M\in L(\Xz,L^2)$, $\F\in L(\Xz,L^2_0)$ with finite dimensional range such that $\tilde\D\in L(\Vnn,L^2_0)$ defined by $\tilde \D \v:=\D\v+\q\cdot\v+\M\v+\F\v$ is bijective.
\end{lemma}
\begin{proof}
First let $\M\v:=-\mean (\q\cdot\v)$ for which it follows that $\D+\q\cdot+\M\in L(V,L^2_0)$.
The new operator acts now on the same spaces as $\D$ and we can perform a perturbation analysis.
Indeed, $\D$ is bijective and $\q\cdot+\M\in L(V,L^2_0)$ is compact from $V$ to $L^2_0$ due to the continuous embedding $V\hookrightarrow \bH^1$ (and because the range of $\M$ is one-dimensional).
Thus $\D+\q\cdot+\M$ is a Fredholm operator with index zero, i.e.\ the range of $\D+\q\cdot+\M$ is closed and $N:=\dim \ker (\D+\q\cdot+\M) = \codim \ran (\D+\q\cdot+\M)<+\infty$.
However, we have no tool at our disposal to ensure that $N=0$ (which would imply the bijectivity of $\D+\q\cdot+\M$).
Thus we perform an additional modification as follows, where we note that the case $N=0$ is included.
We use that $\inner{\div\v,\div\v'}$ is an equivalent scalar product to $\inner{\v,\v'}_{\bH^1}$ on $\Vnn$.
Let $\psi_n, n=1,\dots,N$ be an orthonormal basis with respect to $\inner{\div \v,\div \v'}$ of $\ker(\D+\q\cdot+\M)$, $\phi_n, n=1,\dots,N$ be an orthonormal basis of $\ran(\D+\q\cdot+\M)^\bot$ and set $\F\v:=\sum_{n=1}^N \phi_n \inner{\div\v,\div\psi_n}$.
Thence
\begin{align*}
\tilde\D\in L(\Vnn,L^2_0),\ \tilde \D \v:=\D\v+\q\cdot\v+\M\v+\F\v
\end{align*}
is bijective.
\end{proof}
{N}ote that $\tilde\D$ is also bounded and well-defined on $\IX$, i.e.\ $\tilde\D\in L(\IX,L^2_0)$.
Although the inverse $\tilde\D^{-1}$ will always be considered in the space $L(L^2_0,\Vnn)$.
{For $\u \in \IX$ we construct a topological decomposition mirroring the one in the homogeneous case, in \eqref{eq:Vnn:Wnn}.
    As $\tilde\D$ is bijective on $L(\Vnn,L^2_0)$ we keep $\Vnn$ as in \eqref{eq:Vnn:Wnn} and define
\begin{equation}
\TWnn:=\{\u\in\IX\colon \tilde\D\u=0\}, \label{eq:Wnn:hetero}
\end{equation}
where we use the tilde to indicate the difference to the homogeneous case.
The projections onto $\Vnn$ and $\TWnn$ are now given by
\begin{equation*}
    \widetilde P_V \u:=\tilde\D^{-1}\tilde\D\u, \qquad P_{\TWnn}\u:=\u-\widetilde P_V\u,
\end{equation*}
note that, while $\Vnn$ is the same as in the homogeneous case, the projection $\widetilde P_V $ is different, now defined with respect to $\tilde\D$.
Now $\Vnn\oplus\TWnn$ is again a topological decomposition of $\Xz$.
We keep using the abbreviations $\v:=\widetilde{P_\Vnn}\u$, $\w:=P_{\TWnn}\u$ for $\u\in\Xz$.}

Since $\v\in\Vnn$ it holds
\begin{align*}
\|\div \v\|_{L^2} \geq \beta \|\nabla \v\|_{(L^2)^{3\times3}}.
\end{align*}
Further it follows that
\begin{equations}\label{eq:divq_comp}
(\div+\q\cdot)\w&=(\div+\q\cdot)\u-(\div+\q\cdot)\v\\
&=(\div+\q\cdot)\u-(\div+\q\cdot+\M+\F)\v+(\M+\F)\v\\
&=(\div+\q\cdot)\u-(\div+\q\cdot+\M+\F)\u+(\M+\F)\v\\
&=-(\M+\F)\u+(\M+\F)\v\\
&=-(\M+\F)(\u-\v)\\
&=-(\M+\F)\w
\end{equations}
is a compact operator, which is almost as good as being zero.
Hence the decomposition $\u=\v+\w$ satisfies our wishes.
Thus we build
\begin{align}\label{eq:Tp}
T\u:=\v-\w.
\end{align}
Let $\lambda_-(\m) \in L^\infty$ be the smallest eigenvalue of the symmetric matrix 
\begin{equation}\label{eq:m}
    \m:=-\rho^{-1}\hess(p)+\hess(\phi).
\end{equation}
Further let
\begin{align}\label{eq:theta}
C_M&:=\max\Big\{0, \sup_{x\in\dom} \frac{-\lambda_-(\m(x))}{\gamma(x)}\Big\}
\quad\text{and}\quad
\theta:=\arctan(C_M/|\omega|)\in[0,\pi/2)
\end{align}
for $\omega\neq0$.
\begin{cor}\label{lem:wTc2}
Let $\|c_s^{-1}\bflow\|_{\bL^\infty}^2<\beta^2 \frac{\ull{c_s}^2\ull{\rho}}{\ol{c_s}^2\ol{\rho}} \frac{1}{1+\tan^2\theta}$.
Then $A$ is weakly right $T$-coercive.
\end{cor}
\begin{proof}
Using $T$ as defined in \cref{eq:Tp} we can split $AT=B+K$ with $B,K\in L(\IX)$ given by
\begin{equations}\label{eq:Bp}
\spl B\u,\u' \spr_\IX:=
&-\inner{\rho i\partial_\bflow\v,\opd \w'}+\inner{\rho \opd \w,i\partial_\bflow\v'}\\
&+\inner{\rho\opd \w,\opd \w'}
-\inner{\rho i\partial_\bflow\v,i\partial_\bflow \v'}\\
&+i\omega\inner{\gamma\rho\w,\w'}
+\inner{\rho\m \w,\w'}
+\inner{c_s^2\rho\div \v,\div \v'}\\
& +\inner{c_s^2\rho (\q\cdot \w),(\q\cdot \w')}
+\inner{\rho \F\w,\F \w'}
+\inner{\rho \M\w,\M \w'}
\end{equations}
and
\begin{align*}
    \spl K\u,\u' \spr_\IX:=
&-\inner{\rho (\omega+i\angvel\times) \v, (\omega+i\angvel\times) \v'}
-\inner{\rho (\omega+i\angvel\times) \v, i\partial_\bflow \v'}
-\inner{\rho i\partial_\bflow \v, (\omega+i\angvel\times) \v'}\\
&-\inner{\rho (\omega+i\angvel\times) \v, \opd \w'}
+\inner{\rho \opd \w, (\omega+i\angvel\times) \v'}\\
&-i\omega\inner{\gamma\rho\v,\w'}
+i\omega\inner{\gamma\rho\w,\v'}
-i\omega\inner{\gamma\rho\v,\v'}
+\inner{\rho\m \w,\v'}
-\inner{\rho\m \v,\w'}
-\inner{\rho\m \v,\v'}
\\
&+\inner{c_s^2\rho\q\cdot\v,\div\v'}
+\inner{c_s^2\rho\q\cdot\v,\div\v'}
-\inner{c_s^2\rho\div\w,\q\cdot\v'}
+\inner{c_s^2\rho\q\cdot\v,\div\w'}\\
&-\inner{c_s^2\rho(\div+\q\cdot)\w,(\div+\q\cdot) \w'} 
-\inner{c_s^2\rho(\div+\q\cdot)\w,\div \v'}
+\inner{c_s^2\rho\div\v,(\div+\q\cdot) \w'}\\
&-\inner{\rho \F\w,\F \w'}
-\inner{\rho \M\w,\M \w'}
\end{align*}
for all $\u,\u'\in\IX$.
The operator $K$ is compact due to the compact Sobolev embedding from $V\subset\bH^1$ to $\bL^2$, due to the compactness of $\F,\M$ and \cref{eq:divq_comp}.
Next we show that $B$ is coercive.
Let $\tau\in(0,\pi/2-\theta)$.
First we note that 
\begin{align*}
\frac{1}{\cos(\theta+\tau)}
&\Re \left( e^{-i(\theta+\tau)\sign\omega}
\left( \spl \rho(i\omega\gamma+\m)\w,\w\spr
\right)\right)\\
&=\Re \left(\spl \rho(i\omega\gamma+\m)\w,\w\spr 
\right)
+\sign\omega\tan(\theta+\tau) \Im \left(\spl \rho(i\omega\gamma+\m)\w,\w\spr
\right)\\
&=\spl \rho\m\w,\w\spr
+|\omega| \tan(\theta+\tau) \spl \rho\gamma\w,\w\spr
\geq\spl \rho\lambda_-(\m)\w,\w\spr
+|\omega| \tan(\theta+\tau) \spl \rho\gamma\w,\w\spr\\
&\geq |\omega| (\tan(\theta+\tau)-\tan\theta) \spl \rho\gamma\w,\w\spr,
\end{align*}
whereat the last estimate is due to the definition of $\theta$ \eqref{eq:theta}.
We compute
\begin{align*}
    \frac{1}{\cos(\theta+\tau)}\Re \Big( e^{-i(\theta+\tau)\sign\omega} \spl B\u,\u \spr_\IX\Big)
&\geq\inner{c_s^2\rho\div \v,\div \v}
-\inner{\rho i\partial_\bflow\v,i\partial_\bflow \v}\\
&\quad+\inner{\rho\opd \w,\opd \w}\\
&\quad+(\tan(\theta+\tau)-\tan\theta) |\omega|\inner{\gamma\rho\w,\w}\\
&\quad+\inner{\rho \F\w,\F \w}
+\inner{\rho \M\w,\M \w}
+\inner{c_s^2\rho \q\cdot \w,\q\cdot \w}
\\
&\quad+2\tan(\theta+\tau)\sign\omega \Im\left( \inner{\rho i\partial_\bflow\v,\opd \w} \right).
\end{align*}
We proceed now as in the proof of \Cref{lem:wTc} and estimate
\begin{align*}
    \frac{1}{\cos(\theta+\tau)}\Re \Big( &e^{-i(\theta+\tau)\sign\omega} \spl B\u,\u \spr_\IX\Big)\\
&\geq\inner{c_s^2\rho\div \v,\div \v}
-\big(1+(1-\epsilon)^{-1}\tan^2(\theta+\tau)\big)\inner{\rho i\partial_\bflow\v,i\partial_\bflow \v}\\
&\quad+\epsilon\inner{\rho\opd \w,\opd \w}\\
&\quad+(\tan(\theta+\tau)-\tan\theta) |\omega|\inner{\gamma\rho\w,\w}\\
&\quad+\inner{\rho \F\w,\F \w}
+\inner{\rho \M\w,\M \w}
+\inner{c_s^2\rho \q\cdot \w,\q\cdot \w'}
\end{align*}
The same reasoning as in the proof of \Cref{lem:wTc} yields
\begin{align*}
\inner{c_s^2\rho\div \v,\div \v}
-\big(1+(1-\epsilon)^{-1}\tan^2(\theta+\tau)\big)\inner{\rho i\partial_\bflow\v,i\partial_\bflow \v}
\gtrsim \|\v\|_\IX^2
\end{align*}
and
\begin{align*}
&\epsilon\inner{\rho\opd \w,\opd \w}
+(\tan(\theta+\tau)-\tan\theta) |\omega|\inner{\gamma\rho\w,\w}
\gtrsim \|\partial_\bflow\w\|^2_{\bL^2}+\|\w\|^2_{\bL^2}.
\end{align*}
Using \cref{eq:divq_comp} we know that $\div\w=-\q\cdot\w-\M\w- \F\w$ and we obtain further that
\begin{align*}
&\epsilon\inner{\rho\opd \w,\opd \w}\\
&\quad+(\tan(\theta+\tau)-\tan\theta) |\omega|\inner{\gamma\rho\w,\w} \\
&\quad+\inner{\rho \F\w,\F \w}
+\inner{\rho \M\w,\M \w}
+\inner{c_s^2\rho \q\cdot \w,\q\cdot \w}\\
&\gtrsim
\|\div\w\|^2_{L^2}+\|\partial_\bflow\w\|^2_{\bL^2}+\|\w\|^2_{\bL^2}
=\|\w\|_\IX^2.
\end{align*}
Thus $B$ is uniformly coercive and the proof is finished.
\end{proof}

\subsubsection{The discrete topological decomposition}

Now we mimic this construction on the discrete level.
Let $\PQn$ be the orthogonal projection onto $\Qn$.
Consider the discrete operator
\begin{align*}
\tilde \D_n:=\PQn\tilde \D |_{\Vhnn}.
\end{align*}
Note that $\Vhnn\not\subset\Vnn$ and hence $\tilde \D_n$ is a nonconforming approximation of $\tilde \D$.
{Compared to the homogeneous case we must first ensure that the discrete operator is a suitable approximation of $\tilde \D$.}
\begin{lemma}\label{lem:tDn}
{Let Assumption~\ref{ass:DivStable} be satisfied.} $\tilde \D_n\in L(\Vhnn,\Qn)$ with $p_n=\D_n^{-1}\PQn\div\in L(\Vnn,\Vhnn), \PQn\in L(L^2_0,\Qn)$ forms a discrete approximation scheme of $\tilde\D\in L(\Vnn,L^2_0)$, which approximates $\tilde\D$ and is stable.
In particular, it holds $\lim_{n\to\infty}\|\tilde\D_n^{-1}\PQn\tilde\D\v-\v\|_{\bH^1}=0$ for each $\v\in\Vnn$.
\end{lemma}
\begin{proof}
Due to Lemma~\ref{lem:DnInvPQnDiv} and since $\PQn$ is an orthogonal projection it easily follows that the approximation is a discrete approximation scheme.
For the approximation property we compute for $\v\in\Vnn$
\begin{align*}
\|\tilde\D_n p_n\v-\PQn\tilde\D\v\|_{L^2}
&=\|\tilde\D_n \D_n^{-1}\PQn\div\v-\PQn\tilde\D\v\|_{L^2}\\
&=\|\PQn\tilde\D \D_n^{-1}\PQn\div\v-\PQn\tilde\D\v\|_{L^2}\\
&\leq \|\PQn\tilde\D\|_{L(\bH^1,L^2)} \|\D_n^{-1}\PQn\div\v-\v\|_{\bH^1},
\end{align*}
whereat the right hand-side tends to zero for $n\to\infty$ due to Lemma~\ref{lem:DnInvPQnDiv}.
By construction $\tilde\D\in L(V,L^2_0)$ is bijective and hence the regularity of $\tilde\D_n$ implies its stability.
Since we can split $\tilde\D_n=\D_n+\PQn(q\cdot+\M+\F)|_{\Vhnn}$ into a stable part $\D_n$ and a compact part $\PQn(q\cdot+\M+\F)|_{\Vhnn}$ the regularity of $\tilde\D_n$ follows similarily as in the proof of \Cref{thm:Tcomp}.
The last claim follows from
\begin{align*}
\|\tilde\D_n^{-1}\PQn\tilde\D\v-\v\|_{\bH^1}
\leq \|\tilde\D_n^{-1}\PQn\tilde\D\v - \D_n^{-1}\PQn\div\v\|_{\bH^1}
+\|\D_n^{-1}\PQn\div\v-\v\|_{\bH^1}
\end{align*}
Here the first terms tends to zero, because the discrete approximation scheme of $\tilde\D$ is stable ($p_n=\D_n^{-1}\PQn\div$, $\tilde\D^{-1}\tilde\D\v=\v$), and the second term tends to zero due to Lemma~\ref{lem:DnInvPQnDiv}.
\end{proof}

{Now that we have shown that $\tilde\D_n$ is a suitable approximation of $\tilde\D$ we can proceed similarly as in the homogeneous case, by 
defining a topological decomposition of $\Xn$.
\begin{lemma}
\label{lem:VnWn_hetero}
Let Assumption~\ref{ass:DivStable} be satisfied.
Then the space $\Vhnn$ as in \eqref{eq:Vhnn:Whnn} together with the space 
\begin{equation*}
\TWhnn:=\{\u_n\in \Xn\colon \tilde \D \u_n=0\},
\end{equation*}
form a topological decomposition of $\Xn$ with projections
\begin{align*}
    \widetilde P_{\Vhnn}\u_n:=\tilde\D_n^{-1}P_{Q_n}\tilde\D\u_n, \qquad P_{\TWhnn}\u_n:=\u_n-\widetilde P_{\Vhnn}\u_n,
\end{align*}
being uniformly bounded in $n\in\mathbb{N}$, where for $q_n\in Q_n$ the function $\tilde\D_n^{-1}q_n\in\Vhnn$ is the unique solution to
\begin{align}\label{eq:vn_hetero}
	\text{find }\v_n\in\Vhnn\text{ such that } \tilde\D \v_n = q_n. %
\end{align}
\end{lemma}
\begin{proof}
\cref{lem:tDn} enures that \eqref{eq:vn_hetero} admits a unique solution, and that $\widetilde P_{\Vhnn}$ is uniformly bounded.
Since $\v_n\in\Vhnn$ it also holds
\begin{align*}
\|\div \v_n\|_{L^2}\geq \betad \|\nabla\v_n\|_{(L^2)^{3\times3}}.
\end{align*}
Since $\widetilde P_{\Vhnn}\widetilde P_{\Vhnn}\u_n=\tilde\D_n^{-1}P_{Q_n}\div \widetilde P_{\Vhnn}\u_n=\D_n^{-1}P_{Q_n}P_{Q_n}\div \u_n=\widetilde P_{\Vhnn}\u_n$, $\widetilde P_{\Vhnn}$ is indeed a projection.
Due $\ker \widetilde P_{\Vhnn}=\TWhnn$ the spaces $\Vhnn,\TWhnn$ form indeed a topological decomposition of $\Xn$.
\end{proof}
We abbreviate $\v_n:=\widetilde P_{V_n}\u_n$, $\w_n:=P_{\TWhnn}\u_n$ for $\u_n\in\Xn$.
}
For $\w_n$ we compute
\begin{align}\label{eq:PQn-divq-wn}
\begin{split}
\PQn(\div+\q\cdot)\w_n&=\PQn(\div+\q\cdot)\u_n-\PQn(\div+\q\cdot)\v_n\\
&=\PQn(\div+\q\cdot)\u_n-\PQn(\div+\q\cdot+\M+\F)\v_n+\PQn(\M+\F)\v_n\\
&=\PQn(\div+\q\cdot)\u_n-\PQn(\div+\q\cdot+\M+\F)\u_n+\PQn(\M+\F)\v_n\\
&=-\PQn(\M+\F)\u_n+\PQn(\M+\F)\v_n\\
&=-\PQn(\M+\F)(\u_n-\v_n)\\
&=-\PQn(\M+\F)\w_n,
\end{split}
\end{align}
{which shows that $\PQn(\div+\q\cdot)\w_n$ defines a compact sequence of operators.
This sets up the discrete counterpart of the $T$ operator, and allows us to proceed just like in the homogeneous case, with the following lemma.
\begin{lemma}
\label{lem:Tn_hetero}
Let Assumption~\ref{ass:DivStable} be satisfied.
Then $T_n:=\widetilde P_{\Vhnn}-P_{\TWhnn}=T_n^{-1}\in L(\Xn)$ is uniformly bounded in $n\in\mathbb{N}$.
\end{lemma}
\begin{proof}
Since the spaces $\Vhnn,\TWhnn$ form a topological decomposition of $\Xn$ it follows that $T_n=T_n^{-1}$.
The uniform boundedness of $T_n,T_n^{-1}$ follow from the uniform boundedness of $\widetilde P_{\Vhnn}$, i.e. \cref{lem:tDn}.
\end{proof}
}
\begin{lemma}\label{lem:TnLimit2}
For each $\u\in\IX$ it holds $\lim_{n\to\infty}\|T_n\PXnz\u-\PXnz T\u\|_\IX=0$
\end{lemma}
\begin{proof}
We can proceed similarily as in the proof of Lemma~\ref{lem:TnLimit}.
It suffices to prove $\lim_{n\to\infty}\|\widetilde P_{\Vhnn}\PXnz\u-\PXnz \widetilde P_{\Vnn}\u\|_\IX=0$.
We use $\widetilde P_{\Vhnn}=\tilde\D_n^{-1}\PQn\tilde\D$ and estimate
\begin{align*}
&\|\widetilde P_{\Vhnn}\PXnz\u-\PXnz \widetilde P_{\Vnn}\u\|_\IX\\
&\leq \|\tilde\D_n^{-1}\PQn\tilde\D\PXnz\u-\widetilde P_{\Vnn}\u\|_\IX+\|\widetilde P_{\Vnn}\u-\PXnz \widetilde P_{\Vnn}\u\|_\IX \\
&\leq \|\tilde\D_n^{-1}\PQn\tilde\D\u\!-\!\widetilde P_{\Vnn}\u\|_{\bH^1}
+\|\tilde\D_n^{-1}\PQn\tilde\D\|_{L(\IX,V_n)}\|\u\!-\!\PXnz\u\|_\IX+ \|\widetilde P_{\Vnn}\u\!-\! \PXnz \widetilde P_{\Vnn}\u\|_\IX.
\end{align*}
The second two summands vanish in the limit $n\to\infty$ due to the point-wise convergence of $\PXnz$.
For the first term we apply $\tilde\D\u=\tilde\D \widetilde P_{\Vnn}\u$ and Lemma~\ref{lem:tDn}.
\end{proof}

\subsubsection{Regularity}

Let $\Qnp:=\Qn\oplus\spn\{1\}$ and $\PQnp$ be the $L^2$ orthogonal projection onto $\Qnp$ given by 
$$\PQnp:=\PQn+\M.$$

\begin{lemma}\label{lem:reg2}
If $\|c_s^{-1}\bflow\|_{\bL^\infty}^2<\betad^2 \frac{\ull{c_s}^2\ull{\rho}}{\ol{c_s}^2\ol{\rho}} \frac{1}{1+\tan^2\theta}$, then $(A_n)_{n\in\IN}$ is regular.
\end{lemma}
\begin{proof}
We proceed similarily to the proof of \Cref{lem:reg} and apply \Cref{thm:Tcomp}.
In the previous part of this Section~\ref{subsubsec:RegApr} we already constructed $T_n$ and showed that $T_n\in L(\Xn)$ and $T_n^{-1}=T_n\in L(\Xn)$ are uniformly bounded.
Further, \Cref{lem:TnLimit2} shows that $T_n$ converges pointwise.
Next we split $A_nT_n=B_n+K_n$ into a stable part $B_n\in L(\Xn)$ and a compact part $K_n\in L(\Xn)$.
Recall that $\spl A_nT_n\u_n,\u_n' \spr_\IX=\acow(\v_n-\w_n,\v_n'+\w_n')$.
We start by considering the terms involving $(\div+\q\cdot)$.
Note that
\begin{align*}
\inner{c_s^2\rho\div T_n\u_n,&\div\u_n'}\\
&=\inner{c_s^2\rho\div\v_n,\div\v_n'}
-\inner{c_s^2\rho\div\w_n,\div\v_n'}
+\inner{c_s^2\rho\div\v_n,\div\w_n'}
-\inner{c_s^2\rho\div\w_n,\div\w_n'},
\end{align*}
and
\begin{equation*}
\inner{c_s^2\rho\q\cdot T_n\u_n,\div\u_n'}
=\inner{c_s^2\rho\q\cdot \v_n,\div\v_n'}
\!-\!\inner{c_s^2\rho\q\cdot \w_n,\div\v_n'}
\!+\!\inner{c_s^2\rho\q\cdot \v_n,\div\w_n'}
\!-\!\inner{c_s^2\rho\q\cdot \w_n,\div\w_n'},
\end{equation*}
and
\begin{equation*}
\inner{c_s^2\rho \div T_n\u_n,\q\cdot\u_n'}
=\inner{c_s^2\rho\div \v_n,\q\cdot\v_n'}
\!-\!\inner{c_s^2\rho\div \w_n,\q\cdot\v_n'}
\!+\!\inner{c_s^2\rho\div \v_n,\q\cdot\w_n'}
\!-\!\inner{c_s^2\rho\div \w_n,\q\cdot\w_n'}.
\end{equation*}
Hence
\begin{align}
    &\inner{c_s^2\rho\div T_n\u_n,\div\u_n'}
    +\inner{c_s^2\rho\q\cdot T_n\u_n,\div\u_n'}
    +\inner{c_s^2\rho\div T_n\u_n,\q\cdot\u_n'}\nonumber\\
    &=\inner{c_s^2\rho\q\cdot\v_n,\div\v_n'}
     -\inner{c_s^2\rho\div\w_n,\q\cdot\v_n'}
    +\inner{c_s^2\rho\q\cdot\v_n,\div\w_n'}
    +\inner{c_s^2\rho\div\v_n,\q\cdot\v_n'}\label{eq:divq1}\\
    &\quad-\inner{c_s^2\rho\div\w_n,\div\w_n'}
    -\inner{c_s^2\rho\div\w_n,\q\cdot\w_n'}
    -\inner{c_s^2\rho\q\cdot\w_n,\div\w_n'}\label{eq:divq2}\\
    &
    \quad+\inner{c_s^2\rho\div\v_n,(\div+\q\cdot)\w_n'}
    -\inner{c_s^2\rho(\div+\q\cdot)\w_n,\div\v_n'}
    \label{eq:divq3}\\
    &\quad+\inner{c_s^2\rho\div\v_n,\div\v_n'}\label{eq:divq4}
\end{align}
Line \eqref{eq:divq1} can be moved to the compact operator $K_n$ due to the compact Sobolev embedding from $\Vhnn\subset\bH^1$ to $\bL^2$.
To treat line \eqref{eq:divq2} we note that
\begin{equation*}
\q\cdot\w_n=
\PQn(\q\cdot\w_n)
+\M(\q\cdot\w_n)
+(1-\PQnp)(\q\cdot\w_n)
\end{equation*}
and express
\begin{align*}
\eqref{eq:divq2}&=
-\inner{c_s^2\rho\div\w_n,\div\w_n'}
-\inner{c_s^2\rho\div\w_n,\PQnp(\q\cdot\w_n')}
-\inner{c_s^2\rho\PQnp(\q\cdot\w_n),\div\w_n'}\\
&\quad-\inner{c_s^2\rho\div\w_n,\M(\q\cdot\w_n')}
-\inner{c_s^2\rho\M(\q\cdot\w_n),\div\w_n'}\\
&\quad-\inner{c_s^2\rho\div\w_n,(1-\PQnp)(\q\cdot\w_n')}
-\inner{c_s^2\rho(1-\PQnp)(\q\cdot\w_n),\div\w_n'}\\
&=\inner{c_s^2\rho\PQnp(\q\cdot\w_n),\PQnp(\q\cdot\w_n')}
-\inner{c_s^2\rho(\div+\PQnp\q\cdot)\w_n,(\div+\PQnp\q\cdot)\w_n'}\\
&\quad-\inner{c_s^2\rho\div\w_n,\M(\q\cdot\w_n')}
-\inner{c_s^2\rho\M(\q\cdot\w_n),\div\w_n'}\\
&\quad-\inner{c_s^2\rho\div\w_n,(1-\PQnp)(\q\cdot\w_n')}
-\inner{c_s^2\rho(1-\PQnp)(\q\cdot\w_n),\div\w_n'}\\
&=\inner{c_s^2\rho\PQnp(\q\cdot\w_n),\PQnp(\q\cdot\w_n')}\\
&\quad-\inner{c_s^2\rho \PQnp\F\w_n,\PQnp\F\w_n'}
-\inner{c_s^2\rho\div\w_n,\M(\q\cdot\w_n')}
-\inner{c_s^2\rho\M(\q\cdot\w_n),\div\w_n'}\\
&\quad-\inner{c_s^2\rho\div\w_n,(1-\PQnp)(\q\cdot\w_n')}
-\inner{c_s^2\rho(1-\PQnp)(\q\cdot\w_n),\div\w_n'}
\end{align*}
by means of \eqref{eq:PQn-divq-wn}.
The first line in the former right hand-side is put into $B_n$.
Since $\F$ and $\M$ are compact the second line is put into $K_n$.
The third line tends to zero and is also put into $K_n$.
Indeed, we compute e.g.
\begin{align*}
|\inner{c_s^2\rho(1-\PQnp)(\q\cdot\w_n),\div \v_n'}|
&=|\inner{\q\cdot\w_n,(1-\PQnp)(c_s^2\rho\div \v_n')}|\\
&\leq \norm{\q\cdot\w_n}_{L^2} \norm{(1-\PQnp)(c_s^2\rho\div \v_n')}_{L^2},
\end{align*}
and by means of the discrete commutator property \cite{Bertoluzza99} we estimate
\begin{align}\label{eq:DisCom}
\begin{aligned}
\norm{(1-\PQnp)(c_s^2\rho\div \v_n')}_{L^2}^2
&=\sum_{\tau\in\calT_n} \norm{(1-\PQnp)(c_s^2\rho\div \v_n')}_{L^2(\tau)}^2\\
&=\sum_{\tau\in\calT_n} \norm{(1-\PQnp)((c_s^2\rho-c_\tau)\div \v_n')}_{L^2(\tau)}^2\\
&\leq\sum_{\tau\in\calT_n} \norm{c_s^2\rho-c_\tau}_{L^\infty(\tau)}^2 \norm{\div \v_n'}_{L^2(\tau)}^2\\
&\lesssim h_n^2\sum_{\tau\in\calT_n} \norm{c_s^2\rho}_{W^{1,\infty}(\tau)}^2 \norm{\div \v_n'}_{L^2(\tau)}^2
\end{aligned}
\end{align}
with suitably chosen constants $c_\tau, \tau\in\calT_n$.
Line \eqref{eq:divq3} is treated similarily to line \eqref{eq:divq2}.
Finally, the line \eqref{eq:divq4} is moved to the operator $B_n$.
Hence it holds $A_nT_n=B_n+K_n$ with $B_n$ and $K_n$ defined by
\begin{align*}
\spl K_n\u_n,\u'_n \spr_\IX&:=
-\inner{\rho (\omega+i\angvel\times) \v_n, (\omega+i\angvel\times) \v_n'}
-\inner{\rho (\omega+i\angvel\times) \v_n, i\partial_\bflow \v_n'}
-\inner{\rho i\partial_\bflow \v_n, (\omega+i\angvel\times) \v_n'}\\
&-\inner{\rho (\omega+i\angvel\times) \v_n, \opd \w_n'}
+\inner{\rho \opd \w_n, (\omega+i\angvel\times) \v_n'}\\
&-i\omega\inner{\gamma\rho\v_n,\w_n'}
+i\omega\inner{\gamma\rho\w_n,\v_n'}
-i\omega\inner{\gamma\rho\v_n,\v_n'}\\
&+\inner{\rho\m \w_n,\v_n'}
-\inner{\rho\m \v_n,\w_n'}
-\inner{\rho\m \v_n,\v_n'}\\
&-\inner{\rho \F\w_n,\F \w_n'}
-\inner{\rho \M\w_n,\M \w_n'}\\
\eqref{eq:divq1}\Rightarrow \hspace{1.8em}&\hspace{-1.8em}
\begin{cases}
&+\inner{c_s^2\rho\q\cdot\v_n,\div\v_n'} -\inner{c_s^2\rho\div\w_n,\q\cdot\v_n'} +\inner{c_s^2\rho\q\cdot\v_n,\div\w_n'} +\inner{c_s^2\rho\div\v_n,\q\cdot\v_n'}\\
\end{cases}
\\
\eqref{eq:divq2}\Rightarrow \hspace{1.8em}&\hspace{-1.8em}
\begin{cases}
&-\inner{c_s^2\rho\PQn\F\w_n,\PQn\F \w_n'} \\
&-\inner{c_s^2\rho\div\w_n,(1-\PQnp)(\q\cdot \w_n')}
-\inner{c_s^2\rho(1-\PQnp)(\q\cdot \w_n),\div\w_n'} \\
&-\inner{c_s^2\rho\div\w_n,\M(\q\cdot \w_n')}
-\inner{c_s^2\rho \M(\q\cdot \w_n),\div\w_n'}
\end{cases}
\\
\eqref{eq:divq3}\Rightarrow \hspace{1.8em}&\hspace{-1.8em}
\begin{cases}
&-\inner{c_s^2\rho\PQn\F\w_n,\div \v_n'}
+\inner{c_s^2\rho\div\v_n,\PQn\F \w_n'}\\
&-\inner{c_s^2\rho(1-\PQnp)(\q\cdot\w_n),\div \v_n'}
+\inner{c_s^2\rho\div\v_n,(1-\PQnp)(\q\cdot \w_n')}\\
&-\inner{c_s^2\rho\div\v_n,\M(\q\cdot \w_n')}
-\inner{c_s^2\rho \M(\q\cdot \w_n),\div\v_n'}
\end{cases}
\end{align*}
and
\begin{equations}\label{eq:Bpn}
\spl B_n\u_n,\u'_n \spr_\IX:=
&-\inner{\rho i\partial_\bflow\v_n,\opd \w_n'}+\inner{\rho \opd \w_n,i\partial_\bflow\v_n'}\\
&+\inner{\rho\opd \w_n,\opd \w_n'}
-\inner{\rho i\partial_\bflow\v_n,i\partial_\bflow \v_n'}\\
&+i\omega\inner{\gamma\rho\w_n,\w_n'}
+\inner{\rho\m \w_n,\w_n'}
+\inner{\rho \F\w_n,\F \w_n'}
+\inner{\rho \M\w_n,\M \w_n'}
\\
&+\inner{c_s^2\rho\div \v_n,\div \v_n'}
+\inner{c_s^2\rho \PQn(\q\cdot \w_n),\PQn(\q\cdot \w_n')}
\end{equations}
for all $\u_n,\u_n'\in\Xn$,
{ where $\m$ is as defined in \eqref{eq:m}.}
The operator $(K_n)_{n\in\IN}$ is indeed compact due to the compact Sobolev embedding from $\Vhnn\subset\bH^1$ to $\bL^2$, because $\M, \F$ have a finite dimensional range and because terms involving $1-\PQn$ tend to zero due to \eqref{eq:DisCom}.
It is straightforward to see that $B_n$ is uniformly bounded and that $B_n$ converges pointwise to the operator $B\in L(\IX)$ defined in \cref{eq:Bp}.
It remains to show that $B_n$ is uniformly coercive.
This follows along the lines of the proof of \Cref{lem:wTc2}, whereat $\beta$ is replaced by $\betad$ and we use that $\inner{c_s^2\rho \q\cdot\w_n,\q\cdot\w_n}\lesssim \norm{\w_n}^2_{\bL^2}$.
\begin{extproof}
Let $\tau\in(0,\pi/2-\theta)$.
First we note that 
\begin{align*}
\frac{1}{\cos(\theta+\tau)}
&\Re \left( e^{-i(\theta+\tau)\sign\omega}
\left( \spl \rho(i\omega\gamma+\m)\w_n,\w_n\spr
+\inner{c_s^2\rho \PQn(\q\cdot \w_n),\PQn(\q\cdot \w_n')} 
\right)\right)\\
&=\Re \left(\spl \rho(i\omega\gamma+\m)\w_n,\w_n\spr 
+\inner{c_s^2\rho \PQn(\q\cdot \w_n),\PQn(\q\cdot \w_n')} 
\right)\\
&\quad+\sign\omega\tan(\theta+\tau) \Im \left(\spl \rho(i\omega\gamma+\m)\w_n,\w_n\spr
+\inner{c_s^2\rho \PQn(\q\cdot \w_n),\PQn(\q\cdot \w_n)} 
\right)\\
&=\spl \rho\m\w_n,\w_n\spr
+\inner{c_s^2\rho \PQn(\q\cdot \w_n),\PQn(\q\cdot \w_n')}
+|\omega| \tan(\theta+\tau) \spl \rho\gamma\w_n,\w_n\spr\\
&\geq\spl \rho\lambda_-(\m)\w_n,\w_n\spr
+|\omega| \tan(\theta+\tau) \spl \rho\gamma\w_n,\w_n\spr\\
&\geq |\omega| (\tan(\theta+\tau)-\tan\theta) \spl \rho\gamma\w_n,\w_n\spr,
\end{align*}
whereat the last estimate is due to the definition of $\theta$ \eqref{eq:theta}.
We compute
\begin{align*}
\frac{1}{\cos(\theta+\tau)}\Re \left( e^{-i(\theta+\tau)\sign\omega} \spl B_n\u_n,\u_n \spr_\IX\right)
&\geq\inner{c_s^2\rho\div \v_n,\div \v_n}
-\inner{\rho i\partial_\bflow\v_n,i\partial_\bflow \v_n}\\
&+\inner{\rho\opd \w_n,\opd \w_n}\\
&+(\tan(\theta+\tau)-\tan\theta) |\omega|\inner{\gamma\rho\w_n,\w_n}\\
&+\inner{\rho \F\w_n,\F \w_n}\\
&+2\tan(\theta+\tau)\sign\omega \Im\left( \inner{\rho i\partial_\bflow\v_n,\opd \w_n} \right)
\end{align*}
We proceed now as in the proof of Lemma~\ref{lem:reg} and estimate
\begin{align*}
\frac{1}{\cos(\theta+\tau)}\Re \left( e^{-i(\theta+\tau)\sign\omega} \spl B_n\u_n,\u_n \spr_\IX\right)
&\geq\inner{c_s^2\rho\div \v_n,\div \v_n}
-\big(1+(1-\epsilon)^{-1}\tan^2(\theta+\tau)\big)\inner{\rho i\partial_\bflow\v_n,i\partial_\bflow \v_n}\\
&+\epsilon\inner{\rho\opd \w_n,\opd \w_n}\\
&+(\tan(\theta+\tau)-\tan\theta) |\omega|\inner{\gamma\rho\w_n,\w_n}\\
&+\inner{\rho \F\w_n,\F \w_n}.
\end{align*}
The same reasoning as in the proof of Lemma~\ref{lem:reg} yields
\begin{align*}
\inner{c_s^2\rho\div \v_n,\div \v_n}
-\big(1+(1-\epsilon)^{-1}\tan^2(\theta+\tau)\big)\inner{\rho i\partial_\bflow\v_n,i\partial_\bflow \v_n}
\gtrsim \|\v_n\|_\IX^2
\end{align*}
and
\begin{align*}
\epsilon\inner{\rho\opd \w_n,\opd \w_n}
+(\tan(\theta+\tau)-\tan\theta) |\omega|\inner{\gamma\rho\w_n,\w_n}
\gtrsim \|\partial_\bflow\w_n\|^2_{L^2}+\|\w_n\|^2_{L^2}.
\end{align*}
Since $\div\w_n=-\PQn(\q\cdot\w_n)-\PQn \F\w_n$ we obtain further that
\begin{align*}
\epsilon\inner{\rho\opd \w_n,\opd \w_n}
&+(\tan(\theta+\tau)-\tan\theta) |\omega|\inner{\gamma\rho\w_n,\w_n}
+\inner{\rho \F\w_n,\F \w_n}\\
&\gtrsim
\|\div\w_n\|^2_{L^2}+\|\partial_\bflow\w_n\|^2_{L^2}+\|\w_n\|^2_{L^2}
=\|\w_n\|_\IX^2.
\end{align*}
Thus $B_n$ is uniformly coercive and the proof is finished.
\end{extproof}
\end{proof}

\subsubsection{Convergence}

\begin{theorem}\label{thm:conv2}
Let Assumption~\ref{ass:DivStable} be satisfied and $\|c_s^{-1}\bflow\|_{\bL^\infty}^2<\betad^2 \frac{\ull{c_s}^2\ull{\rho}}{\ol{c_s}^2\ol{\rho}} \frac{1}{1+\tan^2\theta}$.
Let $\u$ be the solution to~\eqref{eq:Galbrun}.
Then there exists an index $n_0>0$ such that for all $n>n_0$ the solution $\u_n$ to \eqref{eq:H1disc} exists
and $\u_n$ converges to $\u$ in the $\IX$-norm with the best approximation estimate
$\|\u-\u_n\|_\IX \lesssim \inf_{\u_n'\in\Xn} \|\u-\u_n'\|_\IX$.
\end{theorem}
\begin{proof}
Due to Lemma~\ref{lem:reg2} the approximation scheme $(A_n)_{n\in\IN}$ is regular.
Since $A$ is bijective~\cite{HallaHohage:21} the claim follows from \Cref{lem:stable,lem:conv}.
\end{proof}
Note that \Cref{rem:convrate} concerning convergence rates still applies.

\begin{remark}\label{rem:th2}
In \cref{rem:th} we already discussed the possibility for different choices of the space $Q_n$ from the one in \eqref{eq:Qn}.
The choice of Taylor-Hood-type discretization using $Q_n = \{f\in L^2_0\colon f|_T\in P_{k-1}(T) \quad\forall T\in\mesh\} \cap H^1$ is also possible in the case of heterogeneous pressure and gravity.
In this case we modify the terms 
\begin{align*}
    \inner{c_s^2\rho(\div + \q\cdot)\u,(\div + \q\cdot)\u'}
-\inner{c_s^2\rho \, \q\otimes \q\u,\u'},
\end{align*}
in \eqref{eq:aq} by inserting the projection $P_{Q_n}$ onto the $H^1$ conforming space of polynomials and obtain
\begin{align*}\inner{c_s^2\rho P_{Q_n}(\div + \q\cdot)\u, P_{Q_n}(\div + \q\cdot)\u'}
    -\inner{c_s^2\rho \,  P_{Q_n}(\q\cdot\u), P_{Q_n}(\q\cdot\u')}.
\end{align*}
As already discussed in \cref{rem:th} this can be implemented using auxiliary variables.
  \end{remark}

\section{Numerical examples}\label{sec:num}
\begin{figure}
    \centering
    \includegraphics[width=.35\textwidth]{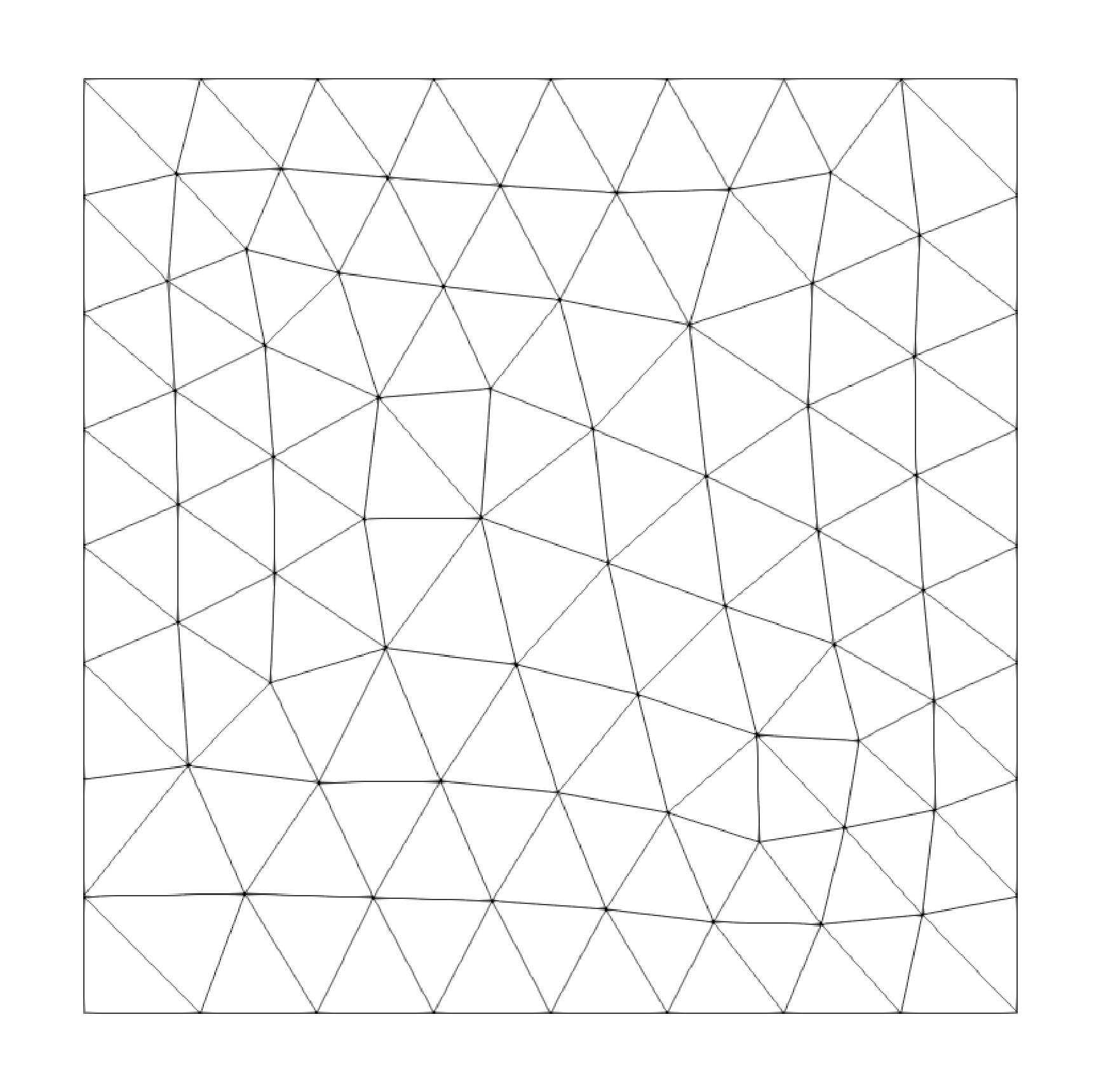}
    \includegraphics[width=.35\textwidth]{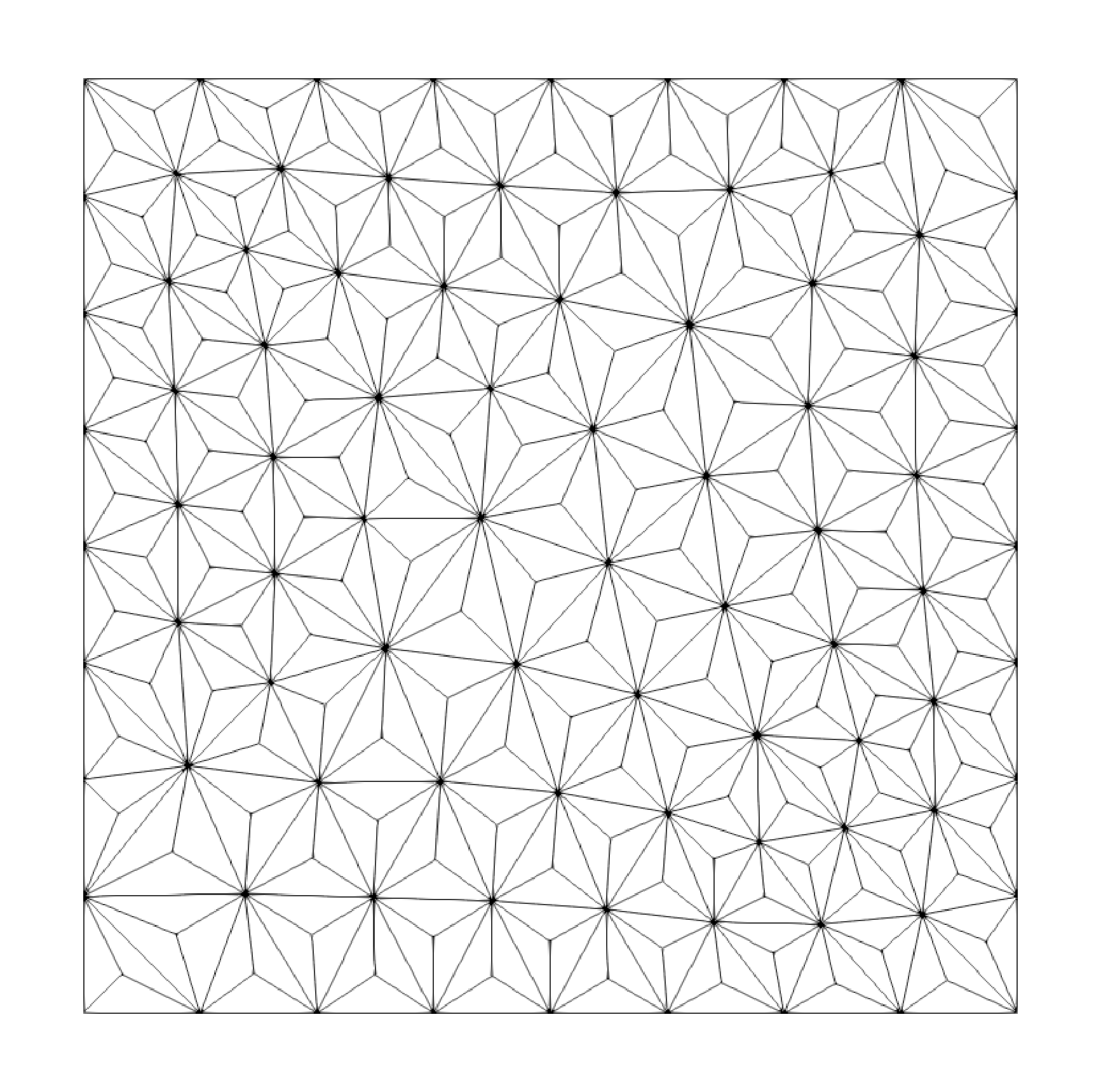}
    \vspace*{-0.3cm}  
    \caption{The coarsest mesh of the sequence of unstructured simplicial meshes is shown on the left, with mesh size $h=0.5$. To construct the second sequence of meshes we apply barycentric mesh refinement, resulting in the mesh on the right.}
    \label{fig:meshes}
\end{figure}

The method has been implemented using \texttt{NGSolve} \cite{ngsolve} and reproduction material is available in \cite{BISNQ9_2023}.
In this section we present numerical examples in the 2D case.
We work with the sesquilinear form given in \cref{eq:cowling} and 
finite element spaces
\begin{align}\label{eq:Xnum}
\Xn:=\{\u\in\bH^1, \u|_\tau \in (P_k(\tau))^2 \quad \forall \tau\in\calT_n\},
\end{align}
with fixed uniform polynomial degree $k\in\IN$.
The error will be measured in the $\norm{\cdot}_\IX$-norm.
We focus on testing the restrictions posed by Assumption~\ref{ass:DivStable} and the smallness assumption on the Mach number $\|c_s^{-1}\bflow\|_{\bL^\infty}$.
In 2D Assumption~\ref{ass:DivStable} requires either: barycentric refinemened meshes and polynomial degree $k\geq 2$ or $k\geq 4$, provided that the meshes have finite degree of degeneracy \cite{ScottVogelius:85}.
To put these conditions to the test, we will consider two sequences of meshes of the domain $\dom=(-4,4)^2$. 
First, shape-regular unstructured simplicial meshes, which include some (nearly-)singular vertices. We will refer to this mesh sequence as unstructured meshes.
These meshes are used to construct the second sequence of meshes. For each mesh in the first sequence we apply barycentric mesh refinement once, constructing the second sequence of meshes.
We will refer to those as the barycentric refined meshes.
A mesh of each type is presented in \Cref{fig:meshes}.
First, we aim to recreate the results obtained in \cite{chabassier:hal-01833043}, which use periodic boundary conditions.
Then, we consider the case of the boundary condition used in this work, given by $\nv\cdot\u=0$.

\subsection{Periodic boundary}
\begin{figure}[!ht]
    \centering
    \hspace*{2em}
    \includegraphics[width=.35\textwidth]{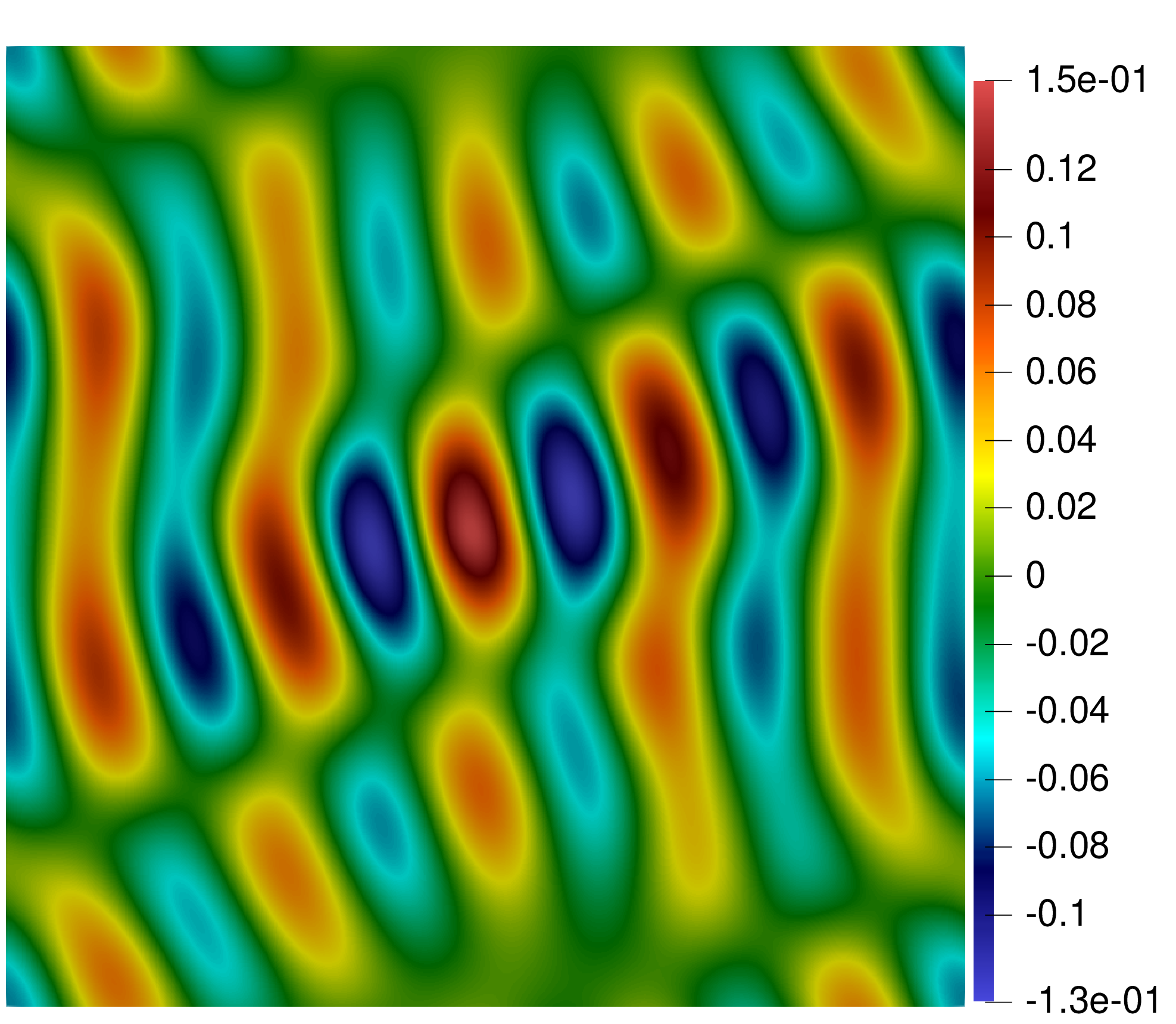}
    \includegraphics[width=.35\textwidth]{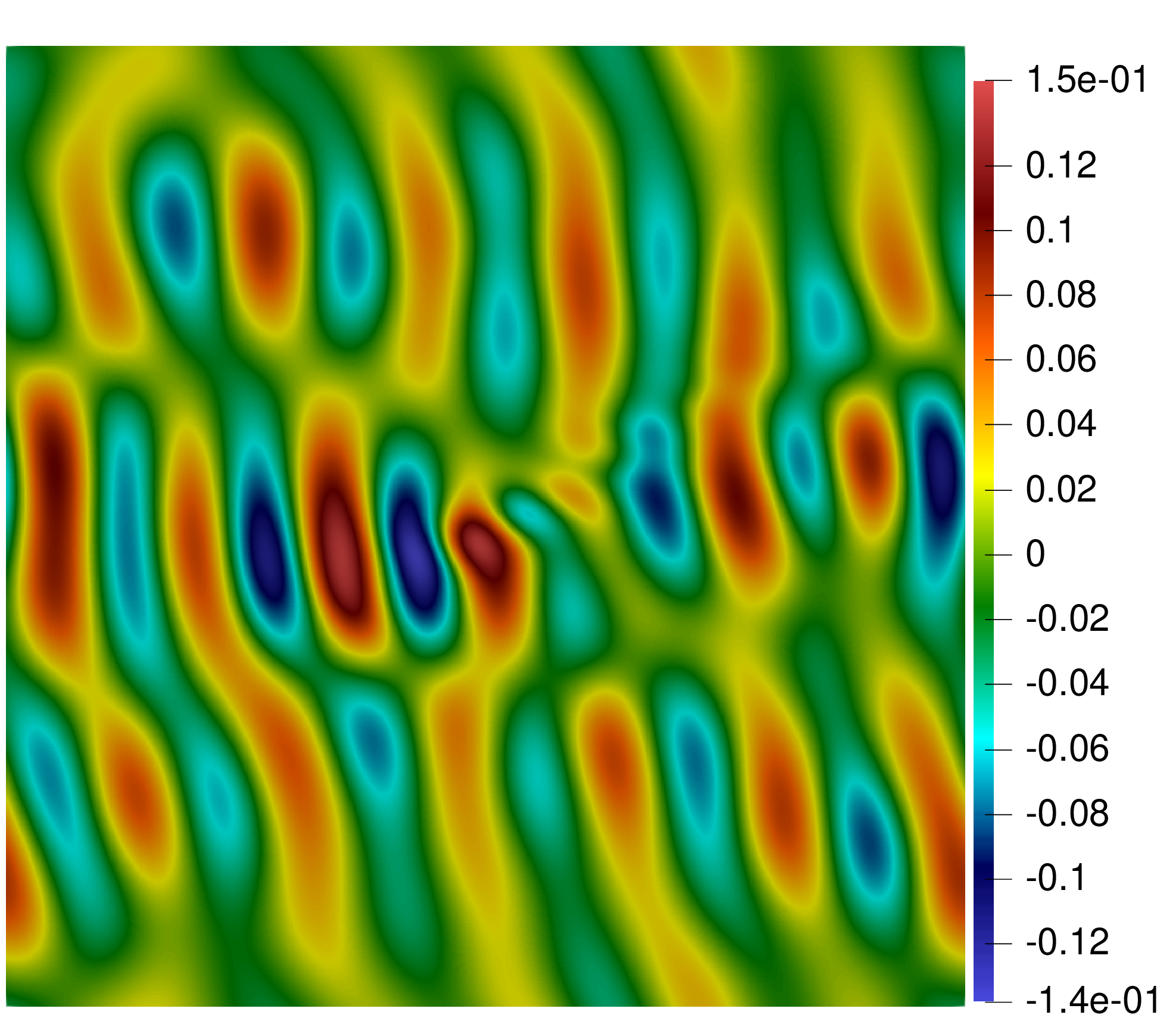}
    \caption{The real part of the first entry of the reference solution computed with $k=5$ and $h=1.5 \cdot 2^{-6}$ for two different values of the coefficient of the flow field $\bflow$, $\coeff=0.2$ on the left and $\coeff=1.5$ on the right.}
    \label{fig:ref_sol}
\end{figure}

We aim to recreate the setting of numerical examples presented in \cite{chabassier:hal-01833043}. 
While we will use the same setting of parameters, there are a few differences.
The major difference is, that our formulation uses a slightly different damping term \cite{HallaHohage:21} than the one considered in  \cite{chabassier:hal-01833043}. 
Furthermore, in \cite{chabassier:hal-01833043} quadrilatera meshes were considered, whereas we will use the simplicial meshes described in \Cref{fig:meshes}.
The setting is as follows:
we consider as computational domain the square $(-4,4)^2$ with periodic boundary conditions, 
and a source term given by 
\begin{align}\label{eq:nums}
    \bff=(-i\omega+\conv)\begin{pmatrix}g\\0\end{pmatrix}
\end{align}
where $g(x,y)$ is the Gaussian given by
$ g(x,y)=\sqrt{a/\pi}\exp(-a(x^2+y^2)).$
Here $a=\log(10^6)$ so that $g$ is equal to $10^{-6}$ on the unit circle.
The parameters are chosen as
\begin{equations}\label{eq:coeffs1}
    &\rho=1.5+0.2\cos(\pi x /4)\sin(\pi y /2),&& c_s^2=1.44+0.16\rho,&& \omega=0.78\times 2\pi, \\ 
    &\gamma=0.1,&& \Omega=(0,0),&& p=1.44\rho + 0.08\rho^2.
\end{equations}
and finally, the background flow is given by
\begin{align}\label{eq:bflow1}
\bflow = \frac{\coeff}{\rho}\begin{pmatrix} 0.3+0.1\cos(\pi y/4)\\ 0.2+0.08\sin(\pi x/4) \end{pmatrix}
\end{align}
The error in the $\norm{\cdot}_\IX$-norm is considered against a reference solution computed with polynomial degree $k=5$ and mesh size $h=1.5\cdot  2^{-6}$.
Plots of the reference solution are shown in \Cref{fig:ref_sol} (compare with \cite[Fig.\ 8, 12]{chabassier:hal-01833043}).

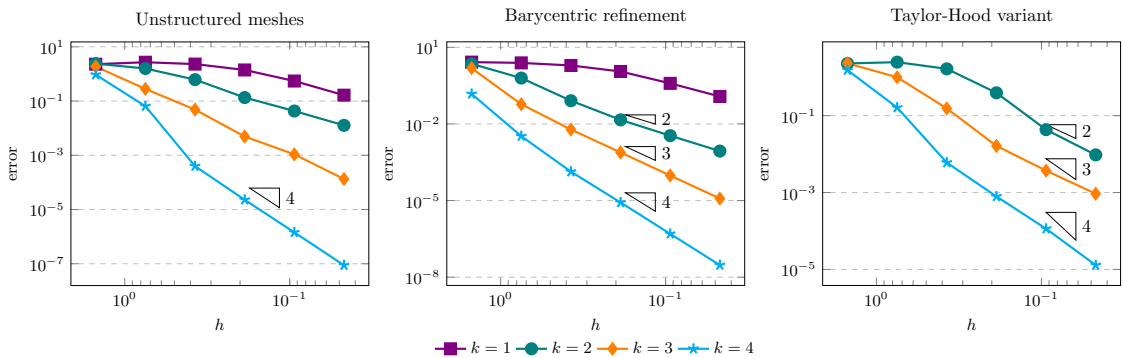
\begin{figure}[ht]
    \centering
    \resizebox{.99\linewidth}{!}{
        \begin{tikzpicture} [spy using outlines={circle, magnification=4, size=2cm, connect spies}]
            \begin{groupplot}[%
                group style={%
                    group name={my plots},
                    group size=3 by 1,
                    horizontal sep=1.3cm,
                },
                legend style={
                    legend columns=6,
                    at={(0.5,-0.2)},
                    anchor=north,
                    draw=none
                },
                ylabel style={at={(-0.18,0.5)}},
                xlabel style={at={(0.5,-0.08)}},
                xlabel={$h$},
                ymajorgrids=true,
                grid style=dashed,
                cycle list name=paulcolors,
                ymin=1e-8, ymax=10,
                ytick={1e-8,1e-6,1e-4,1e-2,1e0},
                ]      
                \nextgroupplot[ymode=log,xmode=log,x dir=reverse, ylabel={error}, title={Unstructured meshes}]
                \foreach \p in {1,2,3,4}{
                    \addplot+[discard if not={order}{\p}] table [x=h, y=error, col sep=comma]{num/ex0.2_per.csv};
                }
                \logLogSlopeTriangle{0.7}{0.1}{0.4}{4}{black}; 
                \nextgroupplot[ymode=log,xmode=log,x dir=reverse, ylabel={error}, title={Barycentric refinement}]
                \foreach \p in {1,2,3,4}{
                    \addplot+[discard if not={order}{\p}] table [x=h, y=error, col sep=comma]{num/ex0.2_bar_per.csv};
                }
                \logLogSlopeTriangle{0.7}{0.1}{0.7}{2}{black}; 
                \logLogSlopeTriangle{0.7}{0.1}{0.57}{3}{black};
                \logLogSlopeTriangle{0.7}{0.1}{0.38}{4}{black}; 
                \legend{$k=1$, $k=2$, $k=3$, $k=4$}
                \nextgroupplot[ymode=log,xmode=log,x dir=reverse, ylabel={error}, title={Taylor-Hood variant}]
                \addplot[draw=none] coordinates {(1,.1)};
                \foreach \p in {2,3,4}{
                    \addplot+[discard if not={order}{\p}] table [x=h, y=error, col sep=comma]{num/ex0.2_per_TH.csv};
                }
                \logLogSlopeTriangle{0.85}{0.1}{0.8}{2}{black}; 
                \logLogSlopeTriangle{0.85}{0.1}{0.65}{3}{black};
                \logLogSlopeTriangle{0.85}{0.1}{0.5}{4}{black}; 
            \end{groupplot}
    \end{tikzpicture}}
    \caption{ Convergence against a reference solution computed with polynomial degree $k=5$ and mesh size $h=1.5\cdot  2^{-6}$.
    We consider the setting described in \eqref{eq:nums} and \eqref{eq:coeffs1} with periodic boundary conditions and
        fixed $\coeff=0.2$ for the background flow $\bflow$ given in \eqref{eq:bflow1}, and different polynomial degree $k=1,2,3,4$ and varying mesh size.
        From left to right we consider: unstructured meshes, barycentric refined meshes, and the Taylor-Hood variant, i.e. unstructured meshes with $Q_n\subset H^1$.
The error is measured in the $\norm{\cdot}_\IX$-norm.
    }
    \label{fig:hconv_per}
\end{figure}

In \Cref{fig:hconv_per} we compare convergence rates for $\alpha=0.1$, putting us safely into the regime of sub-sonic flow.
{Thus satisfying the assumption on the Mach number in \Cref{thm:conv2}. 
If the additional inf-sup stability assumption~\ref{ass:DivStable} is satisfied, then from \cref{thm:conv2} together with \Cref{rem:convrate} we expect convergence rated of order $\calO(h^k)$.
In \Cref{fig:hconv_per} we compare different approaches to satisfy inf-sup stability.
}
We consider the two different types of mesh sequences for polynomial degrees $k=1,2,3,4.$
On unstructured meshes the error is given in \Cref{fig:hconv_per} on the left. There, we observe good convergence rates for $k=4$, after a pre-asymptotic phase, which might be caused by nearly singular vertices. 
For the meshes using barycentric refinement we observe convergence rates of order $\calO(h^{k})$ for $k\geq 2$, shown in \Cref{fig:hconv_per} in the center.
These observations align with the requirements for stability of the Scott-Vogelius element, showing that Assumption~\ref{ass:DivStable} is necessary.
We also show the error for the Taylor-Hood variant which we discussed in \cref{rem:th,rem:th2}, in \Cref{fig:hconv_per} on the right. 
The method used an $H^1$ conforming choice for the space $Q_n$, and we use unstructured meshes.
The method suffers from a long pre-asymptotic phase and shows a worse approximation error compared to the other two methods. 
The rates agree with \Cref{rem:convrate}.

\begin{figure}[ht]
    \centering
    \resizebox{.7\linewidth}{!}{
        \begin{tikzpicture} [spy using outlines={circle, magnification=4, size=2cm, connect spies}]
            \begin{groupplot}[%
                group style={%
                    group name={my plots},
                    group size=2 by 1,
                    horizontal sep=2cm,
                },
                legend style={
                    legend columns=4,
                    at={(-0.2,-0.2)},
                    anchor=north,
                    draw=none
                },
                xlabel={$h$},
                xlabel style={at={(0.5,-0.08)}},
                ymajorgrids=true,
                grid style=dashed,
                cycle list name=paulcolors,
                ]      
                \nextgroupplot[ymode=log,xmode=log,x dir=reverse, ylabel={error}, title={Unstructured meshes}]
                \addplot+[discard if not={order}{4}] table [x=h, y=error, col sep=comma]{num/ex0.2_per.csv};
                \addplot+[discard if not={order}{4}] table [x=h, y=error, col sep=comma]{num/ex0.5_per.csv};
                \addplot+[discard if not={order}{4}] table [x=h, y=error, col sep=comma]{num/ex1.5_per.csv};
                \addplot+[discard if not={order}{4}] table [x=h, y=error, col sep=comma]{num/ex3_per.csv};
                \logLogSlopeTriangle{0.7}{0.1}{0.4}{4}{black}; 
                \nextgroupplot[ymode=log,xmode=log,x dir=reverse, ylabel={consistency error}, title={Unstructured meshes}]
                \addplot+[discard if not={order}{4}] table [x=h, y=conserror, col sep=comma]{num/ex0.2_per.csv};
                \addplot+[discard if not={order}{4}] table [x=h, y=conserror, col sep=comma]{num/ex0.5_per.csv};
                \addplot+[discard if not={order}{4}] table [x=h, y=conserror, col sep=comma]{num/ex1.5_per.csv};
                \addplot+[discard if not={order}{4}] table [x=h, y=conserror, col sep=comma]{num/ex3_per.csv};
                \legend{$\alpha=0.2$,$\alpha=0.5$,$\alpha=1.5$,$\alpha=3$}
            \end{groupplot}
    \end{tikzpicture}}
    \caption{On unstructured meshes with periodic boundary conditions we consider the error in the $\norm{\cdot}_\IX$-norm (left) and consistency-error (right) against a reference solution for different values of the coefficient in the background flow $\bflow$, given in \eqref{eq:bflow1}, and fixed polynomial order $k=4$.}
    \label{fig:hconv_per_coeffs}
\end{figure}
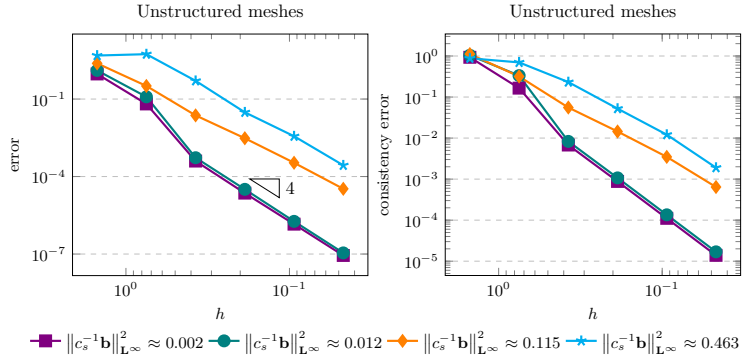

{Now that we observed the importance of using an inf-sup stable methods, we set out to numerically test the smallness requirement on the Mach number in \Cref{thm:conv2}.}
{To do so, we use an inf-sup stable method and we compare different values of the coefficient for the background flow, $\coeff=0.2,0.5,1.5,3$, in \Cref{fig:hconv_per_coeffs}.
We fix $k=4$ and use unstructured meshes.}
We compare against a reference solution computed with $k=5$ and $h=1.5 \cdot 2^{-6}$.
The reference solution for $\alpha=0.2,1.5$ is presented in \Cref{fig:ref_sol}.
As the reference solution changes with $\alpha$, and can give an unreliable comparison, we additionally consider the consistency error, as in \cite{chabassier:hal-01833043}. 
Let us denote 
$$ S_1= 
-\rho\opd^2\u_n
- \gamma \rho i \omega \u_n
,\ \ S_2=
\nabla\left(\rho c_s^2\div \u_n\right)
- (\div \u_n) \nabla p
+\nabla(\nabla p\cdot \u_n) 
- \hess(p)\u_n 
$$
{where the differential operators are applied elementwise.
The consistency error measures the difference between the two terms, which should be (virtually) zero outside the ball $B_{1.5}$ where the source term is located.}
Following \cite[Sect. 3.2]{chabassier:hal-01833043} we define the consistency error by 
$$ \text{consistency error}=\frac{\norm{S_1-S_2}_{\bL^2(\dom\setminus B_{1.5})}}{\norm{S_1}_{\bL^2(\dom\setminus B_{1.5})}}$$
where the error is measured on the domain without the disk with radius $1.5$ centered at the origin, denoted by $\dom\setminus B_{1.5}$. 
Thus removing the effects of the source term $\bff$.

{Estimating the inf-sup constant numerically, see \cref{rem:infsup}, we have $\betad\approx 0.17$ for the considered meshes.
Our assumption on the Mach number in \cref{thm:conv2} then corresponds to
\begin{equation}\label{eq:maex}
    \betad^2 \frac{\ull{c_s}^2\ull{\rho}}{\ol{c_s}^2\ol{\rho}} \frac{1}{1+\tan^2\theta} \approx 0.008.
\end{equation}
The Mach number is approximately $\norm{c_s^{-1}\bflow}_{\bL^\infty}^2\approx 0.002, 0.012, 0.115, 0.463$ for $\coeff=0.2,0.5,1.5,3$.
}
With the choice of $\alpha=1.5$ and $\alpha=3$, we exceed the upper bound notably.
We observe in \Cref{fig:hconv_per_coeffs}, that the error and the consistency error worsen considerably for $\alpha=1.5,3$ and an optimal rate of convergence is not visible for the considered mesh widths.
{On the other hand, for the choice $\alpha=0.5$ we still observe optimal convergence, showing that the bound is not sharp.}
\begin{remark}\label{rem:infsup}
{We have estimated the inf-sup constant $\betad$ in \Cref{ass:DivStable}
 numerically by computing smallest singular value of the matrix 
    $\bM = \bS_q^{-1/2}\bB\bS_{\u}^{-1/2}$, where
    $$(\bS_q)_{i,j}=\inner{q_j,q_i},\qquad 
    (\bB)_{i,j}=\inner{\div\u_j,q_i},\qquad
    \qquad (\bS_{\u})_{i,j}=\inner{\nabla \u_j,\nabla \u_i} + \inner{\u_j,\u_i},$$
    for a basis $(\u_i)_{i=1}^N$ for $\Xn$, as chosen in \eqref{eq:Xnum}, of polynomial degree $k$ and a basis $(q_i)_{i=1}^M$ for $P_{k-1}(\calT)$.
    We recall that the finite element space $\Xn$ does not include any boundary conditions, which is why we chose a stronger norm for $\u$ in the denominator.
}
\end{remark}

\subsection{Normal boundary condition}
In this section we are considering the boundary condition  $\nv\cdot\u=0$.
We do not introduce a new finite element, instead we continue to use the finite element space defined in \eqref{eq:Xnum},
and we incorporate the boundary condition using Nitsche's method. Therefore, we add the following terms to \eqref{eq:cowling}
\begin{align*}
    -\inner{c_s^2 \rho  \u\cdot\nv, \div \u'}_{\partial\dom}
    -\inner{c_s^2\rho  \div \u, \u'\cdot\nv}_{\partial\dom}
    + \frac{\lambda k^2}{h}  \inner{c_s^2 \rho \u \cdot \nv,\u' \cdot \nv}_{\partial\dom} 
\end{align*}
where we choose $\lambda=2^{15}$.
We again consider the domain $\dom=(-4,4)^2$ and the parameters as in \eqref{eq:coeffs1}.
Only the background flow is changed to satisfy $\bflow\cdot\nv=0$ on $\partial\dom$, and will now be given by
\begin{align}\label{eq:bflow2}
\bflow = \frac{\coeff}{\rho}\begin{pmatrix} \sin(\pi x)\cos(\pi y)\\ -\cos(\pi x) \sin(\pi y)  \end{pmatrix}.
\end{align}
The flow additionally fulfills $\div (\rho\bflow)=0$ in $\dom$.
In the following we will consider two examples with different source terms. 

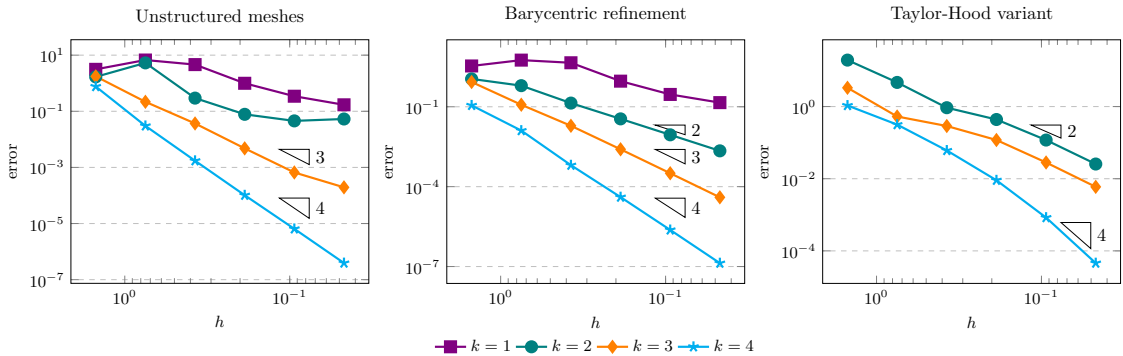
\begin{figure}[ht]
    \centering
    \resizebox{0.99\linewidth}{!}{
        \begin{tikzpicture} [spy using outlines={circle, magnification=4, size=2cm, connect spies}]
            \begin{groupplot}[%
                group style={%
                    group name={my plots},
                    group size=3 by 1,
                    horizontal sep=1.3cm,
                },
                legend style={
                    legend columns=6,
                    at={(0.5,-0.2)},
                    anchor=north,
                    draw=none
                },
                ylabel style={at={(-0.18,0.5)}},
                xlabel style={at={(0.5,-0.08)}},
                xlabel={$h$},
                ymajorgrids=true,
                grid style=dashed,
                cycle list name=paulcolors,
                ymin = 5e-8,
                ymax = 50
                ]      
                \nextgroupplot[ymode=log,xmode=log,x dir=reverse, ylabel={error}, title={Unstructured meshes}]
                \addplot+[discard if not={order}{1}] table [x=h, y=error, col sep=comma] {num/ex0.1_ext.csv};
                \addplot+[discard if not={order}{2}] table [x=h, y=error, col sep=comma] {num/ex0.1_ext.csv};
                \addplot+[discard if not={order}{3}] table [x=h, y=error, col sep=comma] {num/ex0.1_ext.csv};
                \addplot+[discard if not={order}{4}] table [x=h, y=error, col sep=comma] {num/ex0.1_ext.csv};
                \logLogSlopeTriangle{0.8}{0.1}{0.55}{3}{black}; 
                \logLogSlopeTriangle{0.8}{0.1}{0.35}{4}{black}; 
                \nextgroupplot[ymode=log,xmode=log,x dir=reverse, ylabel={error}, title={Barycentric refinement}]
                \addplot+[discard if not={order}{1}] table [x=h, y=error, col sep=comma] {num/ex0.1_ext_bar.csv};
                \addplot+[discard if not={order}{2}] table [x=h, y=error, col sep=comma] {num/ex0.1_ext_bar.csv};
                \addplot+[discard if not={order}{3}] table [x=h, y=error, col sep=comma] {num/ex0.1_ext_bar.csv};
                \addplot+[discard if not={order}{4}] table [x=h, y=error, col sep=comma] {num/ex0.1_ext_bar.csv};
                \logLogSlopeTriangle{0.8}{0.1}{0.65}{2}{black}; 
                \logLogSlopeTriangle{0.8}{0.1}{0.5}{3}{black}; 
                \logLogSlopeTriangle{0.8}{0.1}{0.3}{4}{black}; 
                \legend{$k=1$, $k=2$, $k=3$, $k=4$}
                \nextgroupplot[ymode=log,xmode=log,x dir=reverse, ylabel={error}, title={Taylor-Hood variant}]
                \addplot[draw=none] coordinates {(1,.1)};
                \foreach \p in {2,3,4}{
                    \addplot+[discard if not={order}{\p}] table [x=h, y=error, col sep=comma]{num/ex0.1_ext_TH.csv};
                }
                \logLogSlopeTriangle{0.8}{0.1}{0.78}{2}{black};
                \logLogSlopeTriangle{0.9}{0.1}{0.45}{4}{black}; 
            \end{groupplot}
    \end{tikzpicture}}
    \caption{ Convergence towards an exact solution given in \eqref{eq:ex_sol} for different polynomial orders $k=1,2,3,4$ for varying mesh sizes $h$.
        We consider homogeneous normal boundary contition with fixed $\alpha=0.1$ in the background flow $\bflow$, given in \cref{eq:bflow2}.
        From left to right we consider: unstructured meshes, barycentric refined meshes, and unstructured meshes with $Q_n\subset H^1$.
        The error is measured in the $\norm{\cdot}_\IX$-norm.
    }
    \label{fig:hconv_bnd_ext}
\end{figure}

{As in the periodic study, we start again with a low Mach number flow, that satisfies the assumption in \Cref{thm:conv2}, 
and compare different approaches to satisfy inf-sup stability. }
To this end, we consider convergence against a manufactured solution, by choosing the source term $\bff$ such that the solution will be given by
\begin{align}\label{eq:ex_sol}
    \frac{1}{\rho} \begin{pmatrix}(1+i)g \\ -(1+i)g \end{pmatrix}
\end{align}
where $g$ is again the Gaussian with $a=\log(10^6)$. 
As $g$ equals $10^{-6}$ on the unit circle we can consider the boundary conditions fulfilled numerically to a reasonable degree.
Results for fixed $\alpha=0.1$ are shown in \Cref{fig:hconv_bnd_ext}.
We consider unstructured meshes and barycentric refined meshes. 
Further we include the Taylor-Hood variant outlined in \cref{rem:th,rem:th2} using unstructured meshes and an $H^1$ conforming choice for the space $Q_n$.
For unstructured and barycentric refined meshes we observe the expected convergence rates for $k\geq 3$ and $k\geq 2$, respectively.
The method with $Q_n\subset H^1$ shows again a long pre-asymptotic phase and a worse approximation error compared to the other two methods. Furthermore, for $k=3$ we only observe a long preasymptotic phase, optimal convergence rate is never reached.

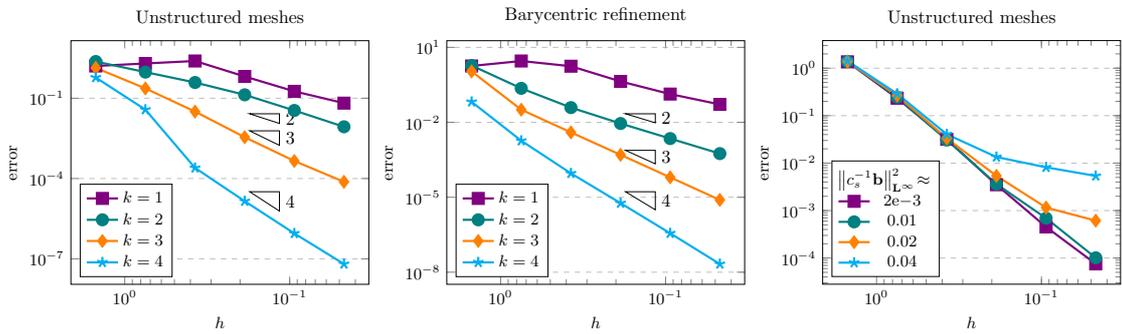
\begin{figure}[ht]
    \centering
    \resizebox{0.99\linewidth}{!}{
        \begin{tikzpicture} [spy using outlines={circle, magnification=4, size=2cm, connect spies}]
            \begin{groupplot}[%
                group style={%
                    group name={my plots},
                    group size=3 by 1,
                    horizontal sep=1.3cm,
                },
                ylabel style={at={(-0.18,0.5)}},
                xlabel style={at={(0.5,-0.08)}},
                xlabel={$h$},
                ymajorgrids=true,
                grid style=dashed,
                cycle list name=paulcolors,
                ]      
                \nextgroupplot[ymin=1e-8,ymax=10,ymode=log,xmode=log,x dir=reverse, ylabel={error}, legend pos=south west, title={Unstructured meshes}]
                \addplot+[discard if not={order}{1}] table [x=h, y=error, col sep=comma] {num/ex0.1.csv};
                \addplot+[discard if not={order}{2}] table [x=h, y=error, col sep=comma] {num/ex0.1.csv};
                \addplot+[discard if not={order}{3}] table [x=h, y=error, col sep=comma] {num/ex0.1.csv};
                \addplot+[discard if not={order}{4}] table [x=h, y=error, col sep=comma] {num/ex0.1.csv};
                \logLogSlopeTriangle{0.7}{0.1}{0.91}{2}{black}; 
                \logLogSlopeTriangle{0.7}{0.1}{0.63}{3}{black};
                \logLogSlopeTriangle{0.7}{0.1}{0.38}{4}{black}; 
                \legend{$k\!\!=\!\!1$, $k\!\!=\!\!2$, $k\!\!=\!\!3$, $k\!\!=\!\!4$}
                \nextgroupplot[ymin=1e-8,ymax=10,ymode=log,xmode=log,x dir=reverse, ylabel={error}, legend pos=south west, title={Barycentric refinement}]
                \addplot+[discard if not={order}{1}] table [x=h, y=error, col sep=comma] {num/ex0.1_bar.csv};
                \addplot+[discard if not={order}{2}] table [x=h, y=error, col sep=comma] {num/ex0.1_bar.csv};
                \addplot+[discard if not={order}{3}] table [x=h, y=error, col sep=comma] {num/ex0.1_bar.csv};
                \addplot+[discard if not={order}{4}] table [x=h, y=error, col sep=comma] {num/ex0.1_bar.csv};
                \logLogSlopeTriangle{0.7}{0.1}{0.7}{2}{black}; 
                \logLogSlopeTriangle{0.7}{0.1}{0.55}{3}{black};
                \logLogSlopeTriangle{0.7}{0.1}{0.38}{4}{black}; 
                \legend{$k\!\!=\!\!1$, $k\!\!=\!\!2$, $k\!\!=\!\!3$, $k\!\!=\!\!4$}
                \nextgroupplot[ymode=log,xmode=log,x dir=reverse, ylabel={error}, legend pos=south west, title={Unstructured meshes}]
                \addplot+[discard if not={order}{3}] table [x=h, y=error, col sep=comma] {num/ex0.1.csv};
                \addplot+[discard if not={order}{3}] table [x=h, y=error, col sep=comma] {num/ex0.2.csv};
                \addplot+[discard if not={order}{3}] table [x=h, y=error, col sep=comma] {num/ex0.3.csv};
                \addplot+[discard if not={order}{3}] table [x=h, y=error, col sep=comma] {num/ex0.4.csv};
                \legend{$\coeff\!=\!0.1$, $\coeff\!=\!0.2$, $\coeff\!=\!0.3$, $\coeff\!=\!0.4$}
            \end{groupplot}
    \end{tikzpicture}}
    \caption{Convergence against a reference solution on the unstructured meshes on the left and barycentric refined meshes in the middle, considering polynomial degree $k=1,2,3,4$. Here $\alpha=0.1$ and the background flow is as in \eqref{eq:bflow2}.
    On the right we consider the unstructured meshes with $k=4$ and different background flows.}
    \label{fig:hconv_bnd}
\end{figure}

Second, we consider again the source term given in \eqref{eq:nums}, this time including the boundary condition, and compare against a reference solution, computed using $k=5$ and $h=1.5\cdot 2^{-6}$ in \Cref{fig:hconv_bnd}.
{Before we investigate the behavior for larger Mach numbers we test convergence against the reference solution for different mesh types and polynomial orders.}
The first two plots in \Cref{fig:hconv_bnd} we fix $\alpha=0.1$ and consider two different mesh types.
For both methods we observe good convergence rates of order $\calO(h^{k})$, however, barycentric refinement show more stable rates and a better error overall.

{Next, we put the assumption on the Mach number in \Cref{thm:conv2} to the test.}
In \Cref{fig:hconv_bnd}, on the right, we consider unstructured meshes, and different values of the coefficient $\alpha$ for for the flow given in \eqref{eq:bflow2}.
We choose $\alpha=0.1,0.2,0.3,0.4$ resulting in the corresponding Mach numbers $\norm{c_s^{-1}\bflow}_{\bL^\infty}^2\hspace{-0.2em}\approx 2\mathrm{e}{-3},\ 0.01,\ 0.02,\ 0.04$.
{From \eqref{eq:maex} we recall that the bound on the Mach number is approximately $0.008$.
True to \Cref{thm:conv2} with the assumptions fulfilled in the case $\norm{c_s^{-1}\bflow}_{\bL^\infty}^2\hspace{-0.2em}\approx 0.01$ we observe the rates given in \Cref{rem:convrate}.
Similar to the periodic case, we still observe convergence for $\norm{c_s^{-1}\bflow}_{\bL^\infty}^2\hspace{-0.2em}\approx 0.01$, even though it is larger than our estimated bound. Nonetheless, for $\norm{c_s^{-1}\bflow}_{\bL^\infty}^2\hspace{-0.2em}\approx 0.02,0.04$ we observe a loss of optimal convergence.}

\section{Conclusion}\label{sec:conclusion}

In this article we reported in \Cref{thm:Tcomp} a new T-compatibility criterion to obtain the regularity of approximations.
As an example of application we considered the damped time-harmonic Galbrun's equation (which is used in asteroseismology) and we proved in \Cref{thm:conv2} convergence for discretizations with divergence stable (Assumption~\ref{ass:DivStable}) $\bH^1$ finite elements.
Although the results of this article constitute only a first step in the numerical analysis for the oscillations of stars.
The subsonic Mach number assumption
\begin{align}\label{eq:subsonic}
 \|c_s^{-1}\bflow\|_{\bL^\infty}^2<\betad^2 \frac{\ull{c_s}^2\ull{\rho}}{\ol{c_s}^2\ol{\rho}}
\end{align}
is far from being optimal.
In stars the density decays with increasing radius and hence the ratio $\frac{\ull{c_s}^2\ull{\rho}}{\ol{c_s}^2\ol{\rho}}$ becomes very small.
Thus a goal is to get rid of this factor in \eqref{eq:subsonic} by a more refined analysis or possibly by more sophisticated discretization methods.
In addition it is desired to replace in \eqref{eq:subsonic} the discrete inf-sup constant of the divergence $\betad$ with a better constant closer to $1$.
The reported computational examples serve only to illustrate the convergence of the finite element method and computational experiments with realistic parameters for stars are eligible.
In particular, a numerical realization of a transparent boundary condition is necessary \cite{Halla:21GalExt,HohageLehrenfeldPreuss:21,BarucqEtal:21}.
Finally we aim to apply the new T-compatibility technique to a number of equations/discretizations for which \cite{Halla:21Tcomp} is too rigid.

\section*{Acknowledgment}
This work was funded by DFG SFB 1456 project 432680300. The first author was supported by DFG project 468728622 and acknowledges that parts of the work was conducted at the Johann Radon Institute for Computational and Applied Mathematics.

\bibliographystyle{abbrvurl}
\bibliography{bib}

\appendix
\section{Construction of \texorpdfstring{$\Hz$}{H1v0}-conforming finite element space}
\label{sec:Xn:construction}
The convenient way to obtain a vectorial $\bH^1$ finite element space is to use a scalar $H^1$ finite element space $Y_n$ and to use $(Y_n)^3$.
Hence if $u_j$ and $\dof_j(u)$ are the basis functions and degrees of freedom of $Y_n$, then $u_j\be_l$ and $\dof_j(\be_l\cdot\bu)$, $l=1,2,3$ with Cartesian unit vectors $\be_j$ are the basis functions and degrees of freedom of $(Y_n)^3$.
However, with this construction it is not clear how to handle the boundary condition $\nv\cdot\bu=0$ and hence the question how to construct finite element spaces of $\Hz$ remains.
To solve this issue for each $j$ we reorganize $u_j\be_l$, $\dof_j(\be_l\cdot\bu)$, $l=1,2,3$ into tangential basis functions and DoFs $\bu_j^{\mathrm{tan},l}$, $\dof_j^\mathrm{tan}$, $l=1,\dots,L_j$ and nontangential ones $\bu_j^{\mathrm{nontan},l}$, $\dof_j^{\mathrm{nontan},l}$, $l=1,\dots,3-L_j$.
Here $L_j=0$ if $\dof_j$ is a vertex DoF associated to a vertrex of $\partial\dom$, $L_j=1$ if $\dof_j$ is a vertex or edge DoF associated to an edge of $\partial\dom$, and $L_j=2$ if $\dof_j$ is a vertex, edge of face DoF associated to a face of $\partial\dom$.
Note that the tangential basis functions $\bu_j^{\mathrm{tan},l}$ will satisfy $\nv\cdot\bu_j^{\mathrm{tan},l}=0$, whereas the nontangential basis functions $\bu_j^{\mathrm{nontan},l}$ will in general satisfy neither $\nv\times\bu_j^{\mathrm{nontan},l}=0$ nor $\nv\cdot\bu_j^{\mathrm{nontan},l}=0$.
However, $\nv\cdot\bu=0$ will imply $\dof_j^{\mathrm{nontan},l}(\bu)=0$.
For $L_j=2$ we simply choose tangential vectors $\bt_1, \bt_2$ and the normal vector $\nv$ and set $\bu_j^{\mathrm{tan},l}=u_j \bt_l$, $\dof_j^{\mathrm{tan},l}(\bu):=\dof_j(\bt_l\cdot\bu)$, $l=1,2$ and $\bu_j^\mathrm{nontan,1}=u_j \nv$, $\dof_j^\mathrm{nontan,1}(\bu):=\dof_j(\nv\cdot\bu)$.
For $L_j=1$ we choose the tangential vector $\bt$ associated to the edge of $\partial\dom$ and $\nv_1, \nv_2$ as the normal vectors of the two adjacent faces.
Thence we set $\bu_j^\mathrm{tan,1}=u_j \bt$, $\dof_j^\mathrm{tan,1}(\bu):=\dof_j(\bt\cdot\bu)$, $l=1,2$ and $\bu_j^{\mathrm{nontan},l}=u_j \nv_l$, $\dof_j^{\mathrm{nontan},l}(\bu):=\dof_j(\nv_l\cdot\bu)$, $l=1,2$.
For $L_j=0$ there exist no tangential DoFs.
Let $\nv_1, \nv_2, \nv_3$ be the normal vectors of the three adjacent faces.
Thence we set $\bu_j^{\mathrm{nontan},l}=u_j \nv_l$, $\dof_j^{\mathrm{nontan},l}(\bu):=\dof_j(\nv_l\cdot\bu)$, $l=1,2,3$.
Thus to obtain a $\Hz$ conforming finite element space we simply set the DoFs associated to the nontangential DoFs to zero.
Hence for $\bu\in\Hz\cap\bH^s$ the obtained finite element space has the same approximation properties as $(Y_n)^3$.

\end{document}